\documentclass[11pt]{amsart}
\usepackage{amsmath, amssymb, latexsym, graphicx, mathrsfs, enumerate, amsthm, textcomp, url, tikz-cd,bbm}
\usepackage[T1]{fontenc}
\usepackage{mathdots}
\usepackage[toc,page]{appendix}
\usepackage{hyperref}
\hypersetup{
    colorlinks=true,
    linkcolor=blue,
    filecolor=magenta,
    urlcolor=cyan,
}

\usepackage{color}

\input xy
\xyoption{all}

\allowdisplaybreaks
\numberwithin{equation}{section}
\newtheorem{theorem}{Theorem}[section]
\newtheorem{proposition}[theorem]{Proposition}
\newtheorem{corollary}[theorem]{Corollary}
\newtheorem{lemma}[theorem]{Lemma}

\newtheorem{definition}[theorem]{Definition}

\newtheorem*{remark*}{Remark}

\newtheorem*{defff}{Definition A}
\newtheorem*{remmm}{Remark B}

\newtheorem{remark}[theorem]{Remark}

\newcommand{\mfo}{\mathfrak{o}}
\newcommand{\mfp}{\mathfrak{p}}

\newcommand{\mfr}{\mathfrak{r}}
\newcommand{\mff}{\mathfrak{f}}

\newcommand{\Mfo}{\mathcal{O}}

\newcommand{\gl}{\mathfrak{gl}}
\newcommand{\cl}{\mathrm{Cl}}
\newcommand{\clb}{\overline{\mathrm{Cl}}}

\newcommand{\End}{\mathrm{End}}

\newcommand{\Z}{\mathbb{Z}}

\newcommand{\ord}{\mathrm{ord}}
\newcommand*{\Rom}[1]{\expandafter\@slowromancap\romannumeral #1@}
\begin{document}

\title[Orbital integrals and ideal class monoids for a Bass order]
{Orbital integrals and ideal class monoids for a Bass order}


\author[Sungmun Cho]{Sungmun Cho}
\author[Jungtaek Hong]{Jungtaek Hong}
\author[Yuchan Lee]{Yuchan Lee}
\thanks{We are supported by  Samsung Science and Technology Foundation under Project Number SSTF-BA2001-04.}

\address{Sungmun Cho \\ Department of Mathematics, POSTECH, Pohang 37673, Republic of Korea}

\email{sungmuncho12@gmail.com}

\address{Jungtaek Hong \\ Department of Mathematics, POSTECH, Pohang 37673, Republic of Korea}

\email{jungtaekhong123@gmail.com}

\address{Yuchan Lee \\ Department of Mathematics, POSTECH, Pohang 37673, Republic of Korea}

\email{yuchanlee329@gmail.com}

\begin{abstract}
A Bass order is an order in a number field for which every fractional ideal can be generated by two elements. Examples include quadratic orders, orders containing the maximal order of a subfield $F$ with $[E:F]=2$, and orders whose discriminant is fourth-power-free.

We prove a closed conductor formula for the number of fractional ideals of a Bass order $R$ modulo multiplication by invertible ideals. For a Bass order, this number is also the number of overorders of $R$, and we give an explicit conductor parametrization of all overorders. The proof combines the classification of local Bass overorders with a local--global argument. Orbital integrals give the corresponding weighted mass formulas, including the split local case. We also prove, by a smoothening procedure, a geometric orbital integral theorem for the relevant integral model; this theorem is presented separately in Section \ref{sec:smoothening}.
\end{abstract}

\maketitle

\tableofcontents

\section{Introduction}
We establish a closed formula, in terms of the conductor, for the number of fractional ideals of a Bass order $R$ modulo multiplication by invertible ideals; see Theorem \ref{thm:intromainthm}. We also give an explicit enumeration of all \emph{overorders} of $R$, that is, orders in $E$ containing $R$.

Some special cases were already known. The quadratic case is classical: every order in a quadratic number field $E$ is uniquely of the form
$R_f=\mathbb Z+f\mathcal O_E$ for some integer $f\geq 1$
\cite[Chapter~20, \S88]{HilbertANT}. The overorders of $R_f$ are exactly
$R_d=\mathbb Z+d\mathcal O_E$ with $d\mid f$; hence, if
$f=\prod_i p_i^{a_i}$, their number is $\prod_i(a_i+1)$.
The present paper proves an analogous conductor formula and parametrization
for Bass orders in arbitrary degree. The subsequent work \cite{AMS} gives
a broader arithmetic classification of overorders; the precise relation is
explained in Remark \ref{rmk:comparisonAMS}.

The proof has three main stages. First, we analyze local Bass overorders
and their conductors (see Theorem \ref{thm:alllocalresults}(2) and
Corollary \ref{cor:sorclbargeneralsplit}). Second, where needed, we use
an exhaustion argument. This combines an independent computation of the
total mass with the orbital integral mass identity and the unit-index
formula (see Propositions \ref{prop:orderidealcounting} and
\ref{rmk:counting}). Finally, a local--global argument yields the global
formula (see Proposition \ref{prop:global_local} and Corollary
\ref{cor:prop_global_local}).

The geometric theorem of Section \ref{sec:smoothening} is also of
separate interest. It adapts the smooth-integral-model method of
\cite{CKL} to the additional strata needed here and computes the
corresponding orbital integrals. We therefore distinguish the arithmetic
classification, the mass identity, and the geometric computation
throughout the paper.

\subsection{Background on a Bass order and an ideal class monoid}
\subsubsection{Bass orders}\label{subsubsec:1.1.1}
A Bass ring was introduced by Hyman Bass in \cite{Bas62} to answer the question:

(\cite[(4.1)]{Lam99})\textit{When is it true that any f.g.\  torsion-free module over a (commutative) Noetherian domain $R$ is isomorphic to a direct sum of ideals?}

Bass found in \cite[Theorem 1.7]{Bas62} that if $R$ is a Noetherian domain whose integral closure $\widetilde{R}$ is a finite $R$-module, then this condition holds if and only if every ideal of $R$ is generated by $2$ elements.
A ring (not necessarily a domain) satisfying the latter condition is called a Bass ring.
Many equivalent characterizations of a Bass ring have been discovered, making it a useful object in various areas.

\begin{proposition}(\cite[pages 96-97]{Lam99})\label{prop:intro}
    Let $R$ be a Noetherian one-dimensional reduced commutative ring.
    Let $\widetilde{R}$ be the normalization of $R$ in its fraction field.
    Suppose that $\widetilde{R}$ is a finitely generated $R$-module.
    Then the following are equivalent.
    \begin{enumerate}
        \item $R$ is a Bass ring, i.e. every ideal of $R$ is generated by $2$ elements.
        \item Every ring $\Mfo$ between $R$ and $\widetilde{R}$ is a Gorenstein ring, that is, every fractional $\Mfo$-ideal $I$ with $(I:I)=\Mfo$ is invertible.
        \item $\widetilde{R}$ is generated by $2$ elements as an $R$-module.
        \item $\widetilde{R}/R$ is a cyclic $R$-module.
    \end{enumerate}
Furthermore, if $R$ is a Bass ring, then every f.g.\ torsion-free $R$-module is isomorphic to a direct sum of ideals.
The converse holds if $R$ is a domain.
\end{proposition}

A Bass order $R$ is an order in a number field $E$ (or in a non-Archimedean local field) satisfying the above condition. Bass orders occur in several natural families (cf. Remark \ref{rmk:globalbass}). If $R$ contains the ring of integers of a subfield $F$ with $[E:F]=2$, then $R$ is Bass; quadratic orders are the basic examples (cf. \cite[Section 2.3]{LW}). Also, if the discriminant of $R$ over $\mathbb{Z}$ is fourth-power-free, then $R$ is Bass (cf. \cite[Theorem 3.6]{Gre82}).
In addition, a cubic order of the form $\mathbb{Z}[x]/(x^3-ax^2+(a-1)x-1)$ for some $a\in \mathbb{Z}$ is a Bass order.
The ideal classes of this order are used to classify similarity classes of Cappell--Shaneson matrices.
They play a key role in the conjecture of Cappell and Shaneson in 4-dimensional manifold theory (cf. \cite[page 44]{AR84}).

\subsubsection{Ideal class monoids}

For an order $R$ in a number field $E$, the \textit{Picard group} $\cl(R)$ is the group of classes of invertible fractional ideals of $R$, and the \textit{ideal class monoid} $\clb(R)$ is the monoid of classes of all fractional ideals of $R$.

Although $\overline{\cl}(R)$ is a fundamental object in number theory, the structure of this monoid is less well understood than that of the Picard group. The local--global study of genera of ideals goes back to Dade--Taussky--Zassenhaus \cite{DTZ}; Proposition \ref{prop:global_local} proves the precise local--global statement needed here. Hofmann and Sircana provide an algorithm to compute orders containing $R$ in \cite{HS}, and Marseglia provides an algorithm to compute $\#\overline{\cl}(R)$ in \cite{Ma20}.

\subsection{Main result}
Before stating the main result of this manuscript, we introduce the conductor ideal $\mff(R)$ to be the largest ideal of $\Mfo_E$ contained in $R$ (cf. Definition A.(4)), where  $\Mfo_E$ is the ring of integers of $E$.
For a Bass order $R$, we express the conductor ideal $\mathfrak{f}(R)$ as an $\Mfo_E$-prime ideal decomposition below (cf. Equation (\ref{descofconductorofR})):
\[    \mathfrak{f}(R)=\left(\mathfrak{p}_1^{2l_1}\cdots \mathfrak{p}_r^{2l_r}\right)\cdot
    \left(\mathfrak{q}_1^{m_1}\cdots \mathfrak{q}_s^{m_s}\right)\cdot \left(\left(\mfr_1\widetilde{\mfr}_1\right)^{n_1}\cdots \left(\mfr_t\widetilde{\mfr}_t\right)^{n_t}\right),
    \]
where $\mathfrak{p}_i$, $\mathfrak{q}_j$, $\mathfrak{r}_k$, and $\widetilde{\mfr}_k$ are distinct prime ideals of $\Mfo_E$ with the following properties:
\begin{itemize}
\item for a maximal ideal $w$ of $R$, let $R_w$ be the $w$-adic completion and $E_w$ its total quotient ring;
\item $\mathfrak p_i$ lies over $w_i$, where $R_{w_i}$ is a domain and $[\kappa_{E_{w_i}}:\kappa_{R_{w_i}}]=1$;
\item $\mathfrak q_j$ lies over $w_j$, where $R_{w_j}$ is a domain and $[\kappa_{E_{w_j}}:\kappa_{R_{w_j}}]=2$;
\item $\mathfrak r_k$ and $\widetilde{\mfr}_k$ lie over $w_k$, where $R_{w_k}$ is not a domain.
\end{itemize}
Here $\kappa_{E_{w_i}}$ is the residue field of the ring of integers in $E_{w_i}$ and $\kappa_{R_{w_i}}$ is the residue field of $R_{w_i}$.
If $R_{w_k}$ is not an integral domain, then there are exactly two maximal ideals  $\mathfrak{r}_k$ and $\widetilde{\mfr}_k$ lying over $w_k$  (cf. Remark \ref{rmk:rirredsplit}).
The main result of this manuscript is the following theorem.

\begin{theorem}(Theorem \ref{thm:bassoverorders})\label{thm:intromainthm}
For a Bass order $R$ of a number field $E$,
\begin{enumerate}
    \item
We have the following equation: $            \#\left(\mathrm{Cl}(R)\backslash\overline{\mathrm{Cl}}(R)\right)=
$
    \begin{multline*}
    \prod_{w\mid \mathfrak{f}(R)}\left( S_{K_w}(R_w)+1 \right)=
    \prod_{p\in \mathcal{P}_R}\prod_{w|p}\left(\frac{S_{p}(R_w)}{d_{R_w}}+1\right)=
    \prod_{i=1}^{r}\left(l_i+1\right)\cdot  \prod_{j=1}^{s}\left(m_j+1\right)\cdot \prod_{k=1}^{t}\left(n_k+1\right).
    \end{multline*}
    Here $d_{R_w}=[\kappa_{R_w}:\mathbb{F}_p]$ and
 $\mathcal{P}_R:=\{p \textit{ a rational prime}\mid p \textit{ divides } \frac{disc(R)}{disc(\Mfo_E)}\}$, where $disc(R)$ is the discriminant of $R$ over $\mathbb{Z}$ and likewise for $disc(\Mfo_E)$.

\item Any order of $E$ containing $R$ is of the form $\langle R, I_{l_i', m_j', n_k'}\rangle$ for a unique ideal $I_{l_i', m_j', n_k'}$ where
\[
I_{l_i', m_j', n_k'}=\left(\mathfrak{p}_1^{2l'_1}\cdots \mathfrak{p}_r^{2l'_r}\right)\cdot
    \left(\mathfrak{q}_1^{m'_1}\cdots \mathfrak{q}_s^{m'_s}\right)\cdot \left(\left(\mfr_1\widetilde{\mfr}_1\right)^{n'_1}\cdots \left(\mfr_t\widetilde{\mfr}_t\right)^{n'_t}\right) ~~~ \textit{ with  }  ~~~
    \left\{
    \begin{array}{l}
         0\leq  l_i'\leq l_i;\\
         0\leq m_j'\leq m_j;\\
         0\leq n_k'\leq n_k.
    \end{array}
    \right.
\]
In this case, $I_{l_i', m_j', n_k'}=\mathfrak{f}(R')$ so that $R'=\langle R, \mff(R')\rangle$.
\end{enumerate}
\end{theorem}
Note that  $\langle R, \mff(R')\rangle$ is  the subring of $\Mfo_E$ generated by elements of $R$ and $\mff(R')$.
In the theorem, $K_w$ is the unramified extension of $\mathbb Q_p$ contained in $E_w$ whose residue field is that of $R_w$, and $\Mfo_{K_w}$ is its ring of integers.  Moreover,
\[
S_{K_w}(R_w)=\operatorname{length}_{\Mfo_{K_w}}
\bigl((\Mfo_E\otimes_RR_w)/R_w\bigr),
\qquad
S_p(R_w):=S_{\mathbb Q_p}(R_w);
\]
see Definition \ref{def:generalserreinv}.
We next explain the local--global reduction and the local exhaustion argument.
\subsubsection{Local-global argument}
A starting point for our work is the following proposition:
\begin{proposition}(Proposition \ref{prop:stratforglobal})\label{prop13}
For an order $R$ of a number field $E$,
we have the formulas:
\[
\#\overline{\mathrm{Cl}}(R)=\sum_{R\subset \mathcal{O}\subset \mathcal{O}_E}\#\mathrm{cl}(\mathcal{O})
~~~~~~~~~\textit{ and } ~~~~~~~~~~~
\mathrm{Cl}(R)\backslash \overline{\mathrm{Cl}}(R)
=\bigsqcup_{R\subset \mathcal{O}\subset \mathcal{O}_E}\overline{\mathrm{cl}(\mathcal{O})},
\]
\[
\textit{where }
\left\{
\begin{array}{l}
     \mathrm{cl}(\mathcal{O}):=\{ [I] \in \overline{\mathrm{Cl}}(\mathcal{O}) \mid (I:I)=\mathcal{O} \}=\{ [I] \in \overline{\mathrm{Cl}}(R) \mid (I:I)=\mathcal{O} \};  \\
     \overline{\mathrm{cl}(\mathcal{O})}:=\mathrm{Cl}(R)\backslash \mathrm{cl}(\mathcal{O}).
\end{array}\right.
\]
Here $\mathcal O$ is an order of $E$ containing $R$, called an \textit{overorder} of $R$.
\end{proposition}
If $R$ is a Bass order, then Proposition \ref{prop:intro}(2) yields that $\#\overline{\mathrm{cl}(\mathcal{O})}=1$ (cf. Remark \ref{rmk:globalbass}(2)) so that
\[  \#(\mathrm{Cl}(R)\backslash \overline{\mathrm{Cl}}(R))=\textit{the number of overorders of $R$ (cf. Proposition \ref{cor:ratioofcl}(1))}.
\]

On the other hand, the set $\mathrm{Cl}(R)\backslash \overline{\mathrm{Cl}}(R)$ satisfies the local-global principle described as follows:

\begin{proposition}(Proposition \ref{prop:global_local})\label{prop:intro14}
For an order $R$ of a number field $E$, the  following map is bijective:
$$    \mathrm{Cl}(R)\backslash \overline{\mathrm{Cl}}(R)\longrightarrow
    \prod\limits_{w: \textit{maximal ideal of $R$}}\overline{\cl}(R_w), ~~~~~~~  \{I\}\mapsto \prod_{w: \textit{maximal ideal of $R$}} [I\otimes_R R_w].   $$
\end{proposition}
This yields the following formulation for $\#\left(\mathrm{Cl}(R)\backslash \overline{\mathrm{Cl}}(R)\right)$ (cf. Corollary  \ref{cor:prop_global_local}(1)):
\[
\#(\mathrm{Cl}(R)\backslash \overline{\mathrm{Cl}}(R))=\#\{\textit{overorders of $R$}\}=\prod\limits_{w:\textit{maximal ideal of $R$}}\#\{\textit{overorders of $R_w$}\}.
\]
Here if $R$ is a Bass order, then so is $R_w$ (cf. Proposition \ref{cor:ratioofcl}(2)).
Therefore, the global problem reduces to the corresponding calculation for the local Bass orders $R_w$.

\subsubsection{Local theory: exhaustion argument using orbital integrals}\label{subsubsec:1.2.2}
Let $E$ be a non-Archimedean local field with ring of integers $\Mfo_E$,
and let $R$ be an order in $E$. We begin with the following reformulation
of an observation of Yun \cite[p.~408]{Yun13}:

\begin{proposition}(Proposition \ref{prop:orderidealcounting})\label{prop:introooverorder}
Let
\[
\left\{
\begin{aligned}
X_R&:=\{\text{fractional $R$-ideals}\},\\
\Lambda_E&:=\pi_E^{\mathbb Z}\subset E^\times,
\end{aligned}
\right.
\qquad
\text{where $E^\times=\Lambda_E\times\Mfo_E^\times$ and
$\Lambda_E$ acts on $X_R$ by multiplication.}
\]

\[
\textit{Then } ~~~~~   \#\overline{\mathrm{Cl}}(R)=\sum_{R\subset \mathcal{O} \subset \mathcal{O}_{E}}\#\mathrm{cl}(\mathcal{O}) ~~~  \textit{   and   } ~~~~
    \#(\Lambda_E \backslash X_R)=\sum_{R\subset \mathcal{O} \subset\mathcal{O}_{E}}\#\mathrm{cl}(\mathcal{O})\cdot \#(\mathcal{O}_{E}^{\times}/ \mathcal{O}^{\times}).
\]
\end{proposition}

Here, $\mathrm{cl}(\Mfo):=\{[I]\in\overline{\mathrm{Cl}}(R)\mid (I:I)=\Mfo\}$. If $R$ is Bass, Proposition \ref{prop:intro}(2), together with Remark \ref{rmk:globalbass}(1), gives $\#\mathrm{cl}(\Mfo)=1$; hence
\begin{equation}\label{eq:introlocal}
         \#(\Lambda_E \backslash X_R)=\sum_{R\subset \mathcal{O} \subset\mathcal{O}_{E}}  \#(\mathcal{O}_{E}^{\times}/ \mathcal{O}^{\times}) \textit{ (cf. Proposition \ref{prop:charofclbarforbass})}.
\end{equation}

In the cases where exhaustion is used, the argument proceeds as follows:
\begin{enumerate}
\renewcommand{\labelenumi}{\textbf{Step (\arabic{enumi})}}
    \item Compute $\#(\Lambda_E\backslash X_R)$.
    \item Construct candidate overorders $\mathcal O$ of $R$ and compute their positive weights $\#(\mathcal O_E^\times/\mathcal O^\times)$.
    \item If their total weight equals $\#(\Lambda_E\backslash X_R)$, positivity implies that no overorder is missing.
\end{enumerate}

\subsubsection{Orbital integrals and smoothening}\label{subsubsec:1.2.3}
Let $\mfo$ be the ring of integers of a subfield $F$ of $E$.
If $R$ is a simple extension of $\mfo$ so that $R\cong \mfo[x]/(\phi(x))$ for an irreducible polynomial $\phi(x)\in \mfo[x]$ of degree $n$, then the left-hand side of Equation (\ref{eq:introlocal}) is the following orbital integral:
$$\textit{$\#(\Lambda_E\backslash X_R)=$ the orbital integral of $\mathfrak{gl}_n(\mfo)$ (cf. Remark \ref{rmk:quotientmeaure}).}$$

Remark \ref{rmk:quotientmeaure} identifies this orbital integral with the volume of the conjugacy class in $\mathfrak{gl}_n(\mfo)$ having characteristic polynomial $\phi(x)$. The conjugacy class is interpreted as the set of $\mfo$-points of an affine scheme over $\mfo$. If this scheme is smooth, its volume is determined by the cardinality of its special fiber. In the present setting the natural model is generally singular, and Section \ref{sec:smoothening} applies a smoothening procedure, following the strategy of \cite{CKL}, when $n$ is odd or when $n$ is even and satisfies the technical condition in \eqref{restriction}.

Marseglia gives a useful criterion for Bass orders in \cite[Proposition 4.6]{Ma24}. Using this criterion, the local argument separates into the following parts:
\begin{enumerate}
\item In the exhaustion cases, the orbital integral in \textbf{Step (1)}
is computed either geometrically in Section \ref{sec:smoothening} or by
the low-degree formulas in \cite{CKL}.
    \item The local algebra constructs candidate overorders and Proposition \ref{rmk:counting} computes their unit-index weights.
    \item The mass identity \eqref{eq:introlocal} proves exhaustion as in \textbf{Step (3)}.
\end{enumerate}
In the remaining cases, the local overorders are classified directly by
conductor calculations, without an orbital integral exhaustion argument;
see Remark \ref{rmk:anotherpfofsec3}. The following theorem summarizes
the local formulas for Bass orders that need not be simple extensions of
$\mfo$ (see Definition \ref{def:orbitalintegralquo}).
\begin{theorem}(Theorem \ref{thm:alllocalresults})\label{thm:introlocalmain}
If $\operatorname{char}(F)=0$ or $\operatorname{char}(F)>n$, every Bass order $R$ in $E$ satisfies
    \[
    \left\{
    \begin{array}{l}
         \#(\Lambda_E\backslash X_R)=q^{S(R)}+[\kappa_E:\kappa_R]\cdot (q^{S(R)-d_R}+q^{S(R)-2d_R}+\cdots+ q^{d_R}+1);\\
         \#\overline{\mathrm{Cl}}(R)=S_K(R)+1.
    \end{array}\right.
    \]
\end{theorem}
Here, $q$ is the cardinality of the residue field of $\mfo$.
In the theorem, $\kappa_E$ and $\kappa_R$ are the residue fields of $\Mfo_E$ and $R$, respectively, and $[\kappa_E:\kappa_R]$ is either $1$ or $2$.  The field $K$ is the unramified extension of $F$ contained in $E$ whose residue field is $\kappa_R$, and $d_R=[K:F]$.  Finally, $S_K(R)$ is the length of $\Mfo_E/R$ as an $\Mfo_K$-module, and $S(R):=S_F(R)$; see Definition \ref{def:generalserreinv}.

\begin{remark}\label{remark:cklconj}

 \cite[Conjecture 1.12]{CKL} formulates a conjecture for the second leading term of orbital integrals for $\mathfrak{gl}_n(\mfo)$ (see also \cite[Section 1.3.2]{CKL}).
 Theorem \ref{thm:introlocalmain} provides evidence for this conjecture. In our setting, the notation of \cite[Conjecture 1.12]{CKL} specializes as follows:
\[
\#B(\gamma)=1, ~~~~  S(\gamma_i)=S(\gamma), ~~~ d_i=1, ~~~~ d_{\gamma_i}=2, ~~~ r_i= [\kappa_E:\kappa_R].
\]
Thus the formula in Theorem \ref{thm:introlocalmain} has the predicted second leading term.
\end{remark}

\subsubsection{Revisiting Theorem \ref{thm:intromainthm}}\label{subsubsec1244}
We now return to Theorem \ref{thm:intromainthm}, with $R$ a Bass order in a number field $E$.
Let $|R|$ be the set of maximal ideals in $R$.
We write (cf. Definition \ref{def:irreddecomsplit})
$$|R|=|R|^{irred}\sqcup |R|^{split} ~~~ \textit{ where  }  ~~~ \left\{
      \begin{array}{l}
|R|^{irred}\subset  \{w\in |R| : \textit{ $R_w$ is an integral domain}\};\\
|R|^{split}\subset  \{w\in |R| : \textit{ $R_w$ is not an integral domain}\}.
      \end{array} \right.$$
In Proposition  \ref{cor:ratioofcl}(2), we prove that $R_w$ is a Bass order and a reduced local ring.
Section \ref{subsubsec:1.2.2} covers the case $w\in |R|^{irred}$.
When $w\in |R|^{split}$, we enumerate all overorders of a Bass order $R_w$ in Section \ref{subsection:formulaforRsplit}, following the same method explained in (1) of Section \ref{subsubsec:1.2.3}.
The formula for $\#\left(\mathrm{Cl}(R)\backslash\overline{\mathrm{Cl}}(R)\right)$ is then obtained from Proposition \ref{prop:intro14}.

For the second assertion, we prove that
$\mff(\langle R, I_{l_i',m_j',n_k'}\rangle)=I_{l_i',m_j',n_k'}$.
The first assertion then shows that these orders exhaust all overorders of $R$.

\begin{remark}[Comparison with Arpin--Marseglia--Springer]
\label{rmk:comparisonAMS}
Independently, Arpin, Marseglia, and Springer obtained a broader
arithmetic classification of overorders \cite{AMS}. More precisely,
\cite[Theorem~4.10]{AMS} places each overorder in a unique
multiplicator ladder associated with a fixed singular maximal ideal,
and \cite[Corollary~4.15]{AMS} combines these ladders to parametrize all
overorders by their conductors and count them. This recovers and
generalizes the classification and enumeration of overorders in
Theorem \ref{thm:intromainthm}. Their framework does not require the
order to be a domain and assumes only the existence of the relevant
multiplicator ladders, rather than the Bass condition at every maximal
ideal. Its proof is considerably shorter and does not use orbital
integrals.

Beyond this overlapping arithmetic classification, the present paper
addresses two further aspects not treated in \cite{AMS}.
First, building on
Yun's orbital integral interpretation \cite{Yun13}, we obtain explicit
local formulas; in the geometric cases, the corresponding orbital
integrals are computed by the smoothening construction of Section
\ref{sec:smoothening}. This construction is globalized in
\cite{CHGlobal}, which, under additional hypotheses, gives an explicit
formula for the $\ell$-adic cohomology of the corresponding compactified
Jacobians.

Second, Proposition \ref{prop:intro14} gives the precise local--global
principle used here for the genus set
$\mathrm{Cl}(R)\backslash\overline{\mathrm{Cl}}(R)$. In
\cite{CHLCubic}, this principle is combined with a classification of
arbitrary local cubic overorders to obtain, for Gorenstein cubic orders,
a closed Euler product and an exact count of genera of integral orbits
in Bhargava's $2\times3\times3$ cubes.

The local--global overorder decomposition proved here and the explicit
cubic overorder and Serre-invariant formulas of \cite{CHLCubic} are
used in \cite[Theorem~1.3]{Lee26} to prove the functional equation of
the completed $L$-function attached to every Gorenstein order in a
cubic number field, thereby supplying the analytic input required in
the Deng--Espinosa argument \cite{DE}. Deng and Espinosa subsequently
proved the identity between Yun's zeta function and their
overorder zeta function by an independent method
\cite[Theorem~1.1 and Corollary~9.1]{DEYun}.
\end{remark}

\textbf{Organization.} Sections \ref{section:icmandobi}-\ref{section:orbbass} develop the local theory and classify Bass overorders in the domain case. Section \ref{sec:smoothening} proves a smooth integral-model theorem and computes the corresponding orbital integrals. Sections \ref{sec:uppbddsimple}-\ref{sec:formulabass} establish the local--global principle, treat the split case, and prove the global conductor formula.

\textbf{Acknowledgments.}
We thank Jungin Lee, Stefano Marseglia, and Evan O'Dorney for helpful comments. We also thank the anonymous referee for a careful reading and detailed suggestions that substantially improved the accuracy, organization, and exposition of the paper.

\part{Local Theory}\label{part1}
In Part 1, we investigate orbital integrals and ideal class monoids for a local field.
Let us start with notations which are used throughout Sections \ref{section:icmandobi}--\ref{sec:smoothening}.

\section*{Notations}\label{part1notations}

\begin{itemize}
\item Let $F$ be a non-Archimedean local field with ring of integers
$\mfo$, uniformizer $\pi$, residue field $\kappa$, and $q=\#\kappa$;
its characteristic is arbitrary unless otherwise stated.

\item
For a finite field extension $F'$ of $F$, we denote by
$\pi_{F'}$  a uniformizer of $F'$,  by  $\Mfo_{F'}$ the ring of integers of $F'$, and  by  $\kappa_{F'}$ the residue field of $F'$.

\item For an element $x\in F'$, $\ord_{F'}(x)$ is the normalized additive valuation satisfying $\ord_{F'}(\pi_{F'})=1$.
If $F'=F$, then we sometimes use $\ord(x)$, instead of $\ord_{F}(x)$.

\item We fix a finite field extension $E$ of $F$ of degree $n$.
Let $e$ be the ramification index of $E/F$ and let  $d=[\kappa_{E}:\kappa]$ so that $n=ed$.

\item An order of $E$ is a subring $\Mfo$ of $E$ such that $\Mfo$ contains $\mfo$ and such that $\Mfo\otimes_{\mfo}F=E$.
Then $\Mfo$ is a local domain.
The maximal ideal of $\Mfo$ is written as $\mathfrak{m}_{\Mfo}$.

\item Let $\kappa_{\Mfo}$ be the residue field of $\Mfo$ and let $d_{\Mfo}=[\kappa_{\Mfo}:\kappa]$.
Then $d_{\Mfo}$ divides the integer $d$.

\item We often  use $R$ to stand for an order of $E$ and $\Mfo$ to stand for an order of $E$ containing $R$.
In this case, $\Mfo$ is called an overorder of $R$.
It is well known that there are finitely many overorders of $R$.

\item For an order $R$ of $E$ and for an ideal $I$ of $\Mfo_E$,
$\langle R, I\rangle$ is  the subring of $\Mfo_E$ generated by elements of $R$ and $I$.

\item We say that an order $R$ is determined by an irreducible polynomial $\phi(x)\in \Mfo_{F'}[x]$ if $R\cong \Mfo_{F'}[x]/(\phi(x))$ as rings.

\item For $a\in A$ or $\psi(x) \in A[x]$ with a flat $\mfo$-algebra $A$, $\overline{a}\in A \otimes_{\mfo} \kappa$ or $\overline{\psi(x)}\in A\otimes_{\mfo} \kappa[x]$ is the reduction of $a$ or $\psi(x)$ modulo $\pi$, respectively.

\item A fractional $\Mfo$-ideal $M$ is a finitely generated $\Mfo$-submodule of $E$ such that $M\otimes_{\mfo}F=E$.
The set of fractional ideals is closed under multiplication and thus forms a monoid.

\item    The ideal quotient $(I:J)$ for two fractional $\Mfo$-ideals $I$ and $J$ is defined to be
$    (I:J)=\{x\in E \mid xJ\subset I\}$.
Then $(I:J)$ is also a fractional $\Mfo$-ideal.

\item A fractional $\Mfo$-ideal $M$  is called invertible if there exists a fractional $\Mfo$-ideal $N$  such that $MN=\Mfo$. If it exists, then it is uniquely characterized by $N=\left(\Mfo:M\right)$.
The set of invertible ideals is closed under multiplication and inversion, and hence forms a group.

\item The ideal class group $\mathrm{Cl}(\Mfo)$ of $\Mfo$ is defined to be the group of equivalence classes of invertible $\Mfo$-ideals up to  multiplication by an element of $E^{\times}$.

\item The ideal class monoid $\overline{\mathrm{Cl}}(\Mfo)$\footnote{Our convention of $\mathrm{Cl}(\Mfo)$ and $\overline{\mathrm{Cl}}(\Mfo)$ follows  \cite{Yun13}, whereas \cite{Ma24} uses $\mathrm{Pic(\Mfo)}$ and $\mathrm{ICM}(\Mfo)$ respectively. }
of $\Mfo$ is defined to be the monoid of equivalence classes of fractional $\Mfo$-ideals up to  multiplication by an element of $E^{\times}$.

\end{itemize}

\section{Ideal class monoids and orbital integrals}\label{section:icmandobi}
This section describes ideal class monoids and orbital integrals using the subsets $\mathrm{cl}(\mathcal{O})$ introduced in Definition \ref{def:clmfo}; see Proposition \ref{prop:orderidealcounting}.
We begin by defining several invariants and the orbital integral for an order of $E$.

\begin{definition}\label{def:invariantofO}
    For an overorder $\Mfo$ of $R$ so that  $R\subset \Mfo\subset \Mfo_{E}$,
\begin{enumerate}
    \item     we define the following  invariants of $\Mfo$:
\[
\left\{
\begin{array}{l}
S(\Mfo):=[\Mfo_{E}:\Mfo], \textit{ the length of $\Mfo_E/\Mfo$ as an $\mfo$-module};\\
u(\Mfo):=\min\{\ord_E(m)\mid m\in \mathfrak{m}_{\Mfo}\} \textit{ so that } \mathfrak{m}_{\Mfo}\Mfo_{E}=\pi_E^{u(\Mfo)}\Mfo_{E}.
\end{array}\right.
\]
$S(\Mfo)$ is called \textit{the Serre invariant} following \cite[Section 2.1]{Yun13}.
Definition \ref{def:generalserreinv} extends this invariant to the global setting.

\item The conductor $\mathfrak{f}(\Mfo)$ of $\Mfo$ is the largest ideal of $\Mfo_E$ which is contained in $\Mfo$. In other words, $\mathfrak{f}(\Mfo)=\{a\in \Mfo_E\mid a\Mfo_E\subset \Mfo\}$.
 Define the integer $f(\Mfo)\in \Z_{\geq 0}$ such that $\mathfrak{f}(\Mfo)=\pi_E^{f(\Mfo)}\Mfo_E \left(\subset \Mfo\right)$.
We sometimes call $f(\Mfo)$ the conductor of $\Mfo$, if it does not cause confusion.
\end{enumerate}
\end{definition}

\begin{remark}\label{rem:invariants3}
We continue to assume that $R\subset \Mfo\subset \Mfo_E$.
    \begin{enumerate}
    \item
If $\Mfo=\Mfo_E$, then $S(\Mfo)=f(\Mfo)=0$ and $u(\Mfo)=1$.

\item If $\Mfo\neq\Mfo_E$, then
the definitions give the following relations:
\[
\left\{
\begin{array}{l}
u(\Mfo)\leq f(\Mfo)\leq  e\cdot S(\Mfo) \leq n\cdot f(\Mfo);\\
1\leq u(\Mfo)\leq u(R)\leq e ~(\textit{since } \pi\in \mathfrak{m}_{\Mfo}), ~~1\leq f(\Mfo)\leq f(R).
\end{array}\right.
\]

\item If $\Mfo \neq \Mfo_E$, then
 $f(\Mfo)$ is the smallest integer and $u(\Mfo)$ is the largest integer satisfying the following inclusions:
\[
\pi_E^{f(\Mfo)}\Mfo_E\subset  \mathfrak{m}_{\Mfo}\subset \pi_E^{u(\Mfo)}\Mfo_{E}.
\]
In other words, $f(\Mfo)$ and $u(\Mfo)$ are optimal to bound $\mathfrak{m}_{\Mfo}$ by $\Mfo_E$-ideals.

\item(\cite[Proposition 2.5]{CKL}) If $R$ is determined by $\phi(x) \in \mfo[x]$, then
$$S(R)=\frac{1}{2}\left(\ord\bigl(\operatorname{disc}(\phi)\bigr)-\ord\bigl(\operatorname{disc}(E/F)\bigr)\right).$$
\end{enumerate}
\end{remark}

\begin{definition}\label{def:orbitalintegralquo}
We define the orbital integral for an order $R$ as follows:
\begin{enumerate}
    \item
The orbital integral  for  $\phi(x) \in \mfo[x]$ is
$\#(\Lambda_E \backslash X_R)$,
where
\[
\left\{
\begin{array}{l}
R\cong\mfo[x]/(\phi(x));\\
X_R: \textit{ the set of fractional $R$-ideals};\\
\Lambda_E \left(\subset E^{\times}\right):  \textit{ a free abelian group such that } \Lambda_E= \pi_E^{\mathbb{Z}}.
\end{array} \right.
\]
Here, $\Lambda_E$ is complementary to $\Mfo_{E}^{\times}$ inside $E^{\times}$ and acts on $X_R$ by multiplication.

\item Extending the above, the orbital integral for an order $R$ is defined to be
$\#(\Lambda_E \backslash X_R)$.
\end{enumerate}

\end{definition}

\begin{remark}\label{rmk:quotientmeaure}
\begin{enumerate}
    \item
The orbital integral for $\phi(x)$ is indeed defined to be the integral \[
\mathcal{SO}_{\gamma, d\mu}=\int_{\mathrm{T}_{\gamma}(F)\backslash \mathrm{GL}_n(F)} \mathbf{1}_{\mathfrak{gl}_n(\mathfrak{o})}(g^{-1}\gamma g) d\mu(g),
\]
where $\gamma$ is a matrix in $\mathfrak{gl}_n(\mathfrak{o})$ whose characteristic polynomial is $\phi(x)$,  $T_{\gamma}$ is the centralizer of $\gamma$ in $\mathrm{GL}_n$, and $\mathbf{1}_{\mathfrak{gl}_n(\mathfrak{o})}$ is the characteristic function of $\mathfrak{gl}_n(\mathfrak{o})\subset \mathfrak{gl}_n(F)$.
The measure $d\mu(g)$ is called \textit{the quotient measure}, which is given in \cite[Section 1.3]{Yun13}.
 \cite[Section 2.2.1]{CKL} also gives a self-contained explanation for $d\mu(g)$.

Using the quotient measure $d\mu(g)$, it is well known that  $\mathcal{SO}_{\gamma, d\mu}=\#(\Lambda_E \backslash X_R)$ (cf. \cite[Theorem 2.5 and Corollary 4.6]{Yun13}), which justifies Definition \ref{def:orbitalintegralquo}.

\item \cite[Section 3]{FLN} introduces  another measure, called \textit{the geometric measure}, which is described precisely in Section \ref{subsec:geom}.
Section \ref{sec:smoothening} uses this measure.

    \end{enumerate}

\end{remark}

The description of $\overline{\mathrm{Cl}}(R)$ and $\#(\Lambda_E \backslash X_R)$ in Proposition \ref{prop:orderidealcounting} is based on a subset $\mathrm{cl}(\mathcal{O})$ of $\overline{\mathrm{Cl}}(R)$ defined in Definition \ref{def:clmfo}.
For a Bass order, $\#\mathrm{cl}(\Mfo)=1$ for every overorder $\Mfo$; see Remark \ref{rmk:bassideal1}(1).
Thus the description in terms of $\mathrm{cl}(\Mfo)$ is particularly simple for Bass orders.

\begin{definition}\label{def:clmfo}
    For an overorder $\Mfo$ of $R$ so that  $R\subset \Mfo\subset \Mfo_{E}$,
    we define the set $\mathrm{cl}(\mathcal{O})$ as follows:
    \[\mathrm{cl}(\mathcal{O}):=\{ [I]\in \overline{\mathrm{Cl}}(R)\mid (I:I)=\mathcal{O}\}.\]
\end{definition}
Definition \ref{def:globalclcl} extends this notion to orders in finite \'{e}tale algebras.
If $I'=aI$ for $a\in E^\times$, then $(I':I')=(I:I)$. Thus $(I:I)$ depends only on $[I]$, and $\mathrm{cl}(\mathcal{O})$ is well-defined.
Since $(I:I)$ is the largest order over which $I$ is a fractional ideal, the set $\mathrm{cl}(\mathcal{O})$ also admits the following description:
\[
\mathrm{cl}(\mathcal{O})=\{[I]\in\overline{\mathrm{Cl}}(R)\mid
\mathcal{O}\text{ is the largest order over which }I\text{ is a fractional ideal}\}.
\]

\begin{proposition}\label{prop:orderidealcounting}
We have the following equations:
\[        \#\overline{\mathrm{Cl}}(R)=\sum_{R\subset \mathcal{O} \subset \mathcal{O}_{E}}\#\mathrm{cl}(\mathcal{O}) ~~~  \textit{   and   } ~~~~
    \#(\Lambda_E \backslash X_R)=\sum_{R\subset \mathcal{O} \subset\mathcal{O}_{E}}\#\mathrm{cl}(\mathcal{O})\cdot \#(\mathcal{O}_{E}^{\times}/ \mathcal{O}^{\times}).
\]
\end{proposition}

\begin{proof}
The first equation follows from the disjoint union
\[        \overline{\mathrm{Cl}}(R)=\bigsqcup_{R\subset \mathcal{O} \subset \mathcal{O}_{E}}\{ [I]\in \overline{\mathrm{Cl}}(R)\mid (I:I)=\mathcal{O} \}=\bigsqcup_{R\subset \mathcal{O} \subset \mathcal{O}_{E}}\mathrm{cl}(\mathcal{O}).
\]
For the second, Yun's identity \cite[p.~408]{Yun13} gives
$\#(\Lambda_E\backslash X_R)
=\sum_{[I]\in\overline{\mathrm{Cl}}(R)}
\#\bigl(\Mfo_E^\times/\operatorname{Aut}(I)\bigr)$, where
$\operatorname{Aut}(I):=\{x\in E^\times:xI=I\}$.
Since $\operatorname{Aut}(I)=(I:I)^\times$, grouping the ideal
classes by their multiplicator rings gives the second equation.
\end{proof}

\begin{proposition}[Unit-index formula]\label{rmk:counting}
For every order $\Mfo$ of $E$, we have
\begin{equation}\label{eq:generalunitindex}
\#(\Mfo_E^{\times}/\Mfo^{\times})
=\frac{q^d-1}{q^{d_{\Mfo}}-1}\,
q^{S(\Mfo)-d+d_{\Mfo}}.
\end{equation}
In particular, this index is $q^{S(\Mfo)}$ if
$\kappa_{\Mfo}=\kappa_E$, and is
$(q^{d_{\Mfo}}+1)q^{S(\Mfo)-d_{\Mfo}}$ if
$[\kappa_E:\kappa_{\Mfo}]=2$.
\end{proposition}

\begin{proof}
The formula is immediate for $\Mfo=\Mfo_E$. Suppose that
$\Mfo\ne\Mfo_E$. Since
$\mathfrak m_{\Mfo}=\Mfo\cap\pi_E\Mfo_E$, the exact sequence
$0\to\pi_E\Mfo_E/\mathfrak m_{\Mfo}\to
\Mfo_E/\Mfo\to\kappa_E/\kappa_{\Mfo}\to0$ gives
$\#(\pi_E\Mfo_E/\mathfrak m_{\Mfo})
=q^{S(\Mfo)-d+d_{\Mfo}}$. Reduction on units gives
$\#(\Mfo_E^\times/\Mfo^\times)
=\frac{q^d-1}{q^{d_{\Mfo}}-1}
\#((1+\pi_E\Mfo_E)/(1+\mathfrak m_{\Mfo}))$.

Since $f(\Mfo)\geq1$ and
$\pi_E^{f(\Mfo)}\Mfo_E\subset\mathfrak m_{\Mfo}$, the maps
$1+x\mapsto x$ induce bijections of sets
\[
\frac{1+\pi_E\Mfo_E}
     {1+\pi_E^{f(\Mfo)}\Mfo_E}
\longleftrightarrow
\frac{\pi_E\Mfo_E}
     {\pi_E^{f(\Mfo)}\Mfo_E},
\qquad
\frac{1+\mathfrak m_{\Mfo}}
     {1+\pi_E^{f(\Mfo)}\Mfo_E}
\longleftrightarrow
\frac{\mathfrak m_{\Mfo}}
     {\pi_E^{f(\Mfo)}\Mfo_E}.
\]
Taking the ratio of the resulting cardinality equalities and using the
preceding formulas proves \eqref{eq:generalunitindex}. The two
special cases follow from $d_{\Mfo}=d$ and $d=2d_{\Mfo}$, respectively.
\end{proof}

\section{Formula for orbital integrals and ideal class monoids in a Bass order}\label{section:orbbass}

In this section, we provide explicit formulas for orbital integrals and ideal class monoids of a Bass order.
We also enumerate all overorders of a Bass order using the conductor.

\subsection{Characterization of a Bass order}\label{subsec:charbass}
\begin{definition}\label{def:bass}
For an order $R$ of $E$,
\begin{enumerate}
    \item(\cite[Proposition 3.4]{Ma24}) $R$ is Gorenstein if every fractional $R$-ideal $I$ with $R=(I:I)$ is invertible.

   \item(\cite[Proposition 4.6]{Ma24} or \cite[Theorem 2.1]{LW}) $R$ is called a Bass order if every overorder of $R$ is Gorenstein, equivalently if every ideal is generated by two elements.

    \end{enumerate}
\end{definition}

\begin{remark}\label{rmk:bassideal1}
\begin{enumerate}
    \item
Since $\Mfo$ is a one-dimensional Noetherian local domain,
a fractional $\Mfo$-ideal is
invertible if and only if it is principal by \cite[Proposition~12.4]{Neu}.
Therefore $\Mfo$ is Gorenstein if and only if $\#\mathrm{cl}(\Mfo)=1$.
This yields the following description of a Bass order:
\[     \textit{$R$ is a Bass order if and only if $\#\mathrm{cl}(\Mfo)=1$ for every overorder $\Mfo$ of $R$.}
\]

\item The maximal order $\Mfo_E$ is a Bass order, and
$\#(\Lambda_E\backslash X_{\Mfo_E})
=\#\overline{\mathrm{Cl}}(\Mfo_E)
=\#\mathrm{Cl}(\Mfo_E)=1$.
Thus, in this section, we work with a Bass order $R\ne\Mfo_E$.

\item If $[E:F]=2$, then any order of $E$ is Bass by \cite[Section 2.3]{LW}.
Note that any order in this case is a simple extension of $\mfo$.
We use this fact in Section \ref{section:orbbass}.

\item If $R$ is a Bass order, then any overorder $\Mfo$ of $R$ is  Bass as well, since an overorder of $\Mfo$ is also an overorder of $R$ which is Gorenstein.
\end{enumerate}
\end{remark}

\begin{proposition}\label{prop:charofclbarforbass}
    For a Bass order  $R$,
$$\#\overline{\mathrm{Cl}}(R)=\textit{the number of overorders of } R \textit{  and  }
\#(\Lambda_E \backslash X_R)=\sum_{R\subset \mathcal{O} \subset\mathcal{O}_{E}}  \#(\mathcal{O}_{E}^{\times}/ \mathcal{O}^{\times}).
$$
\end{proposition}
\begin{proof}
    This directly follows from Proposition \ref{prop:orderidealcounting} since $\#\mathrm{cl}(\Mfo)=1$ by Remark \ref{rmk:bassideal1}.
\end{proof}

For examples of Bass orders, see \cite[Section~2.3]{LW} and Remark~\ref{rmk:globalbass}; we now fix notation.

\begin{definition}\label{def:fieldK}
    \begin{enumerate}
    \item    Let $K$ be the unramified extension of $F$ contained in $E$ whose residue field is $\kappa_R$.
Then  $d_R=[\kappa_R\left(=\kappa_{K}\right):\kappa]=[K:F]$.
    By \cite[Lemma 3.1]{CKL}, we have $\Mfo_{K}\subset R$.

\item Define $S_K(\Mfo)$ to be  the length of $\Mfo_{E}/\Mfo$ as an $\Mfo_K$-module.
Then $S_K(\Mfo)=S(\Mfo)/d_R$.

\item Let $n_R=n/d_R=[E:K]$.
\end{enumerate}
\end{definition}

\begin{proposition}(\cite[Corollary 4.4 and Proposition 4.6]{Ma24})\label{prop:bassclass}
For $R\neq \Mfo_E$,         $R$ is a Bass order if and only if either $u(R)=2$ and $\kappa_E=\kappa_R$, or $u(R)=1$  and $[\kappa_E:\kappa_R]=2$.
\end{proposition}

\begin{proof}
    This is a simple application of \cite[Corollary 4.4 and Proposition 4.6]{Ma24}.
    By loc. cit., $R$ is a Bass order if and only if  $\dim_{\kappa_R}\Mfo_E/\pi_E^{u(R)}\Mfo_E=2$. The claim follows from this.
\end{proof}

Therefore a nonmaximal Bass order falls into two cases.
Sections \ref{subsec:1stcase}-\ref{subsec:2ndcase} treat them separately.
For a simple extension, the following proposition gives a polynomial
criterion for being Bass.

\begin{proposition} \label{prop:polynomialofbass}
For $R\neq \Mfo_E$, suppose that $R=\mfo[\alpha]$ is a simple extension of $\mfo$ determined by $\phi(x)\in\mfo[x]$ of degree $n$. Let $g(x)\in K[x]$ be the monic minimal polynomial of $\alpha$ over $K$. Then
\begin{itemize}
\item $g(x)\in \Mfo_K[x]$ has degree $n_R$, and $\overline{g(x)}=(x-\bar a)^{n_R}$ for some $\bar a\in\kappa_R$;
\item for every lift $a\in\Mfo_K$ of $\bar a$, the order $R$ is Bass if and only if either $n_R=2$ or $\ord_K(g(a))=2$.
\end{itemize}
\end{proposition}

\begin{proof}
Since $\Mfo_K\subset R$ and $\kappa_R=\kappa_K$, the element
$\alpha$ is integral over $\Mfo_K$, so $g\in\Mfo_K[x]$, and its
reduction has the stated form. By \cite[Lemma~I.4]{Ser},
$\mathfrak m_R=(\pi,\alpha-a)$; hence
$u(R)=\min\{e,\ord_E(\alpha-a)\}$. Put
$f_{E/K}=[\kappa_E:\kappa_K]$. Since
$g(a)=(-1)^{n_R}N_{E/K}(\alpha-a)$ and $n_R=ef_{E/K}$,
the normalized valuations give
$\ord_K(g(a))=f_{E/K}\ord_E(\alpha-a)$.

Since $R\ne\Mfo_E$, we have $n_R>1$. If $n_R=2$, then $R$ is Bass
by Remark \ref{rmk:bassideal1}(3). Assume that $n_R>2$. If $R$ is
Bass, Proposition \ref{prop:bassclass} gives
$(f_{E/K},u(R))=(1,2)$ or $(2,1)$. In either case $e>u(R)$, so
$\ord_E(\alpha-a)=u(R)$, and the norm identity gives
$\ord_K(g(a))=2$.

Conversely, suppose that $\ord_K(g(a))=2$. Since both factors in the
norm identity are positive integers, we have
$(f_{E/K},\ord_E(\alpha-a))=(1,2)$ or $(2,1)$. Since
$n_R=ef_{E/K}>2$, we have $e>2$ in the first case and $e>1$ in the
second; hence $u(R)=2$ or $1$, respectively. Proposition
\ref{prop:bassclass} therefore shows that $R$ is Bass.
\end{proof}

For every Bass order $R$, since $R$ is Gorenstein, the standard
index--discriminant formula and \cite[Proposition~4]{DCD} give
$\#(\Mfo_E/R)^2=N_{E/K}\bigl(\mathfrak f(R)\bigr)$. Since
$\#(\Mfo_E/R)=q^{d_RS_K(R)}$ and
$N_{E/K}\bigl(\mathfrak f(R)\bigr)
=q^{d_R[\kappa_E:\kappa_R]f(R)}$, we obtain
\begin{equation}\label{eq:relsercon}
f(R)=\frac{2S_K(R)}{[\kappa_E:\kappa_R]}.
\end{equation}

We now state the main theorem of Part \ref{part1}, which gives formulas for $\#(\Lambda_E\backslash X_R)$ and $\#\overline{\mathrm{Cl}}(R)$.
\begin{theorem}\label{thm:alllocalresults}
Assume that $\operatorname{char}(F)=0$ or $\operatorname{char}(F)>n$.
\begin{enumerate}
    \item
    For a Bass order $R$ of $E$, we have
    \[
    \left\{
    \begin{array}{l}
         \#(\Lambda_E\backslash X_R)=q^{S(R)}+[\kappa_E:\kappa_R]\cdot (q^{S(R)-d_R}+q^{S(R)-2d_R}+\cdots+ q^{d_R}+1);\\
         \#\overline{\mathrm{Cl}}(R)=S_K(R)+1=[\kappa_E:\kappa_R]\cdot f(R)/2+1.
    \end{array}\right.
    \]
    Here $K$ is as in Definition \ref{def:fieldK}; if $S(R)=0$, the sum in parentheses is empty.

\item If $R'$ is an overorder of $R$, then $R'=\langle R, \pi_E^{f(R')}\Mfo_E\rangle$. We refer to \nameref{part1notations}  for $\langle R, \pi_E^{f(R')}\Mfo_E\rangle$.

\end{enumerate}
\end{theorem}

\begin{proof}
If $R=\Mfo_E$, both assertions are immediate. Otherwise, Proposition
\ref{prop:bassclass} reduces the proof to the two cases in Sections
\ref{subsec:1stcase} and \ref{subsec:2ndcase}; Equation
\eqref{eq:relsercon} gives the last equality in~(1), while Theorems
\ref{thm:computeorbital39} and \ref{thmorb2nd} and Corollaries
\ref{cor:1stmain} and \ref{cor:3rdmain} give the remaining formulas
in~(1) and assertion~(2).
\end{proof}

\subsection{\texorpdfstring{The case $u(R)=2$ and $\kappa_E=\kappa_R$}
{The case u(R)=2 and the residue fields agree}}
\label{subsec:1stcase}
We treat the first case of Proposition \ref{prop:bassclass} and derive
closed formulas for $\#(\Lambda_E\backslash X_R)$ and
$\#\overline{\mathrm{Cl}}(R)$. Assume throughout that
$R\ne\Mfo_E$, $u(R)=2$, and $\kappa_E=\kappa_R$.
Then $E/K$ is totally ramified of degree $e=n_R>1$, and $K/F$ is
unramified of degree $d=d_R$.

\begin{proposition}\label{prop:bassfirstcase}
\begin{enumerate}
\item $R\supset \Mfo_K[u\pi_E^2]$ for some $u\in\Mfo_E^\times$.
\item If $e=2$, then every order of $E$ containing $\Mfo_K$ is Bass.
\item If $e>2$, then $R$ is a simple extension of $\Mfo_K$ if and only if
$R=\Mfo_K[u\pi_E^2]$ for some $u\in\Mfo_E^\times$.
\end{enumerate}
\end{proposition}

\begin{proof}
Since $u(R)=2$, there exists $u\in\Mfo_E^\times$ such that
$u\pi_E^2\in\mathfrak m_R$; as $\Mfo_K\subset R$, this proves~(1).
Part~(2) follows from Remark \ref{rmk:bassideal1}(3), with $K$ as base
field, and part~(3) from Proposition \ref{prop:polynomialofbass}.
\end{proof}

\subsubsection{The cases of odd $e>2$ and $e=2$}\label{sec:first_oddor2}
\begin{proposition}\label{prop:bassfirstcase_2}
Suppose that $e>2$ is odd or $e=2$.
    \begin{enumerate}
    \item If $e>2$ is odd, then $S(R)\leq d(e-1)/2$, with equality if $R$ is a simple extension of $\Mfo_K$.
\item If $e=2$, then $R=\Mfo_K[u\pi_E^t]$ for some
$u\in\Mfo_E^\times$ and an odd integer $t$, with
$S(R)=d(t-1)/2$.
    \end{enumerate}
\end{proposition}

\begin{proof}
For~(1), Proposition \ref{prop:bassfirstcase}(1) gives an order
$R_0:=\Mfo_K[x]\subset R$, where $x=u\pi_E^2$ for some
$u\in\Mfo_E^\times$; when $R$ is simple, Proposition
\ref{prop:bassfirstcase}(3) shows that $R$ itself is of this form.
Since $e$ is odd, the normalized $E$-valuations of
$1,x,\ldots,x^{\frac{e-1}{2}},
x^{\frac{e+1}{2}}/\pi,\ldots,x^{e-1}/\pi$
run through $0,\ldots,e-1$, so these elements form an
$\Mfo_K$-basis of $\Mfo_E$. Comparing this basis with the power basis
$1,x,\ldots,x^{e-1}$ of $R_0$ gives
$S_K(R_0)=(e-1)/2$. Thus
$S(R)\leq S(R_0)=d(e-1)/2$, with equality when $R$ is simple.

For~(2), since $\Mfo_E=\Mfo_K\oplus\Mfo_K\pi_E$ and $\Mfo_K$ is a
DVR, one has
$R=\Mfo_K\oplus a\Mfo_K\pi_E=\Mfo_K[a\pi_E]$
for some nonzero $a\in\Mfo_K$. Since
$\ord_E(a\pi_E)=2\ord_K(a)+1$, write
$a\pi_E=u\pi_E^t$ for some $u\in\Mfo_E^\times$ and an odd integer
$t$. The element
$u\pi_E^t/\pi^{\frac{t-1}{2}}$ is a uniformizer of $E$, so
$\{1,u\pi_E^t/\pi^{\frac{t-1}{2}}\}$ and
$\{1,u\pi_E^t\}$ are $\Mfo_K$-bases of $\Mfo_E$ and $R$,
respectively. Hence $S_K(R)=(t-1)/2$ and
$S(R)=d(t-1)/2$.
\end{proof}
Note that $S(R)$ is not bounded when $e=2$ since any order in this case is Bass by Remark \ref{rmk:bassideal1}(3).

\begin{theorem}\label{thmorb1st}
Assume that $\operatorname{char}(F)=0$ or
$\operatorname{char}(F)>n$, and suppose that either $e>2$ is odd or
$e=2$. If $R$ is a simple extension of $\Mfo_K$, then
\[
\#(\Lambda_E\backslash X_R)=1+q^d+\cdots+q^{S(R)}.
\]
\end{theorem}

\begin{proof}
For odd $e>3$, the formula follows from Theorem \ref{thmorbital} and
\cite[Lemma~3.2]{CKL}; the cases $e=3$ and $e=2$ follow from
\cite[Remarks~6.7.(1) and~5.7]{CKL}, respectively.
\end{proof}

\begin{theorem}\label{thm:overorder1}
Assume that $\operatorname{char}(F)=0$ or
$\operatorname{char}(F)>n$, and suppose that either $e>2$ is odd or
$e=2$. Let $R$ be a simple extension of $\Mfo_K$. Thus
$R=\Mfo_K[u\pi_E^2]$ if $e>2$ is odd, while
$R=\Mfo_K[u\pi_E^t]$ for an odd integer $t$ if $e=2$, where
$u\in\Mfo_E^\times$; see Propositions
\ref{prop:bassfirstcase}(3) and \ref{prop:bassfirstcase_2}(2).
Set
\[
\Mfo_k:=
\left\{
\begin{array}{ll}
\Mfo_K[u\pi_E^2,\frac{(u\pi_E^2)^k}{\pi}],
& \frac{e+1}{2}\leq k\leq e,\quad e>2,\\
\Mfo_K[\frac{u\pi_E^t}{\pi^{(t-k)/2}}],
& 1\leq k\leq t,\quad k\text{ odd},\quad e=2.
\end{array}
\right.
\]

\begin{enumerate}
\item The orders $\Mfo_k$ are pairwise distinct and exhaust all
overorders of $R$.
\item
$\#\overline{\mathrm{Cl}}(R)
=S_K(R)+1
=\frac{S(R)}{d_R}+1
=
\left\{
\begin{array}{ll}
\frac{e+1}{2},& e>2,\\
\frac{t+1}{2},& e=2.
\end{array}
\right.$
\end{enumerate}
\end{theorem}

\begin{proof}
By Proposition \ref{prop:charofclbarforbass},
$\#(\Lambda_E\backslash X_R)
=\sum_{R\subseteq\Mfo\subseteq\Mfo_E}
\#(\Mfo_E^\times/\Mfo^\times)$.
Each $\Mfo_k$ is an overorder of $R$, and a direct comparison of
$\Mfo_K$-bases gives
$S_K(\Mfo_k)=k-(e+1)/2$ when $e>2$ is odd and
$S_K(\Mfo_k)=(k-1)/2$ when $e=2$. Since
$\kappa_{\Mfo_k}=\kappa_E$, Proposition \ref{rmk:counting} gives
\begin{equation}\label{eq:dlbare22}
\#(\Mfo_E^\times/\Mfo_k^\times)
=
\begin{cases}
q^{d(k-(e+1)/2)},& e>2\text{ odd},\\
q^{d(k-1)/2},& e=2.
\end{cases}
\end{equation}

The candidate weights therefore sum to
$1+q^d+\cdots+q^{d(e-1)/2}$ and
$1+q^d+\cdots+q^{d(t-1)/2}$, respectively. By Theorem
\ref{thmorb1st}, each sum equals the orbital integral
$\#(\Lambda_E\backslash X_R)$ in the corresponding case.
Since every
overorder contributes a positive weight, no additional overorder can
occur. The displayed values of $S_K(\Mfo_k)$ show that the candidates
are pairwise distinct. This proves~(1). Part~(2) follows from the candidate count,
Proposition \ref{prop:charofclbarforbass}, and the displayed values of
$S_K(\Mfo_k)$.
\end{proof}

We now remove the simple-extension hypothesis.

\begin{corollary}\label{cor:1stmain}
Assume that $\operatorname{char}(F)=0$ or $\operatorname{char}(F)>n$,
and that either $e=2$ or $e>2$ is odd.
\[
\left\{
\begin{array}{l}
\#(\Lambda_E\backslash X_R)
=q^{S(R)}+q^{S(R)-d}+\cdots+q^d+1;\\
\#\overline{\mathrm{Cl}}(R)
=S_K(R)+1=\frac{S(R)}{d_R}+1.
\end{array}
\right.
\]
Every overorder $R'$ of $R$ satisfies
$R'=\langle R,\pi_E^{f(R')}\Mfo_E\rangle$.
\end{corollary}

\begin{proof}
If $e>2$ is odd, Proposition \ref{prop:bassfirstcase}(1) gives a
simple order $R_0:=\Mfo_K[u\pi_E^2]\subset R$, which is Bass by
Proposition \ref{prop:bassclass}. If $e=2$, let $R_0:=R$; it is simple
by Proposition \ref{prop:bassfirstcase_2}(2). Thus, in either case,
$R_0$ is a simple Bass order contained in $R$. Theorem
\ref{thm:overorder1} shows that $R=\Mfo_k$ for some $k$ and that its
overorders are precisely the $\Mfo_{k'}$ with $k'\leq k$.

As $k'$ ranges over these overorders, $S_K(\Mfo_{k'})$ ranges over
$0,\ldots,S_K(R)$. Thus \eqref{eq:dlbare22} gives the weights
$1,q^d,\ldots,q^{S(R)}$, and Proposition
\ref{prop:charofclbarforbass} gives the two displayed formulas.

The preceding values of $S_K(\Mfo_{k'})$ and
\eqref{eq:relsercon} show that these overorders are determined by their
conductors. For $R'=\Mfo_{k'}$, set
$A:=\langle R,\pi_E^{f(R')}\Mfo_E\rangle$. Since
$A\subseteq R'$ and $\pi_E^{f(R')}\Mfo_E\subseteq A$, the two
inclusions give $f(A)\geq f(R')$ and $f(A)\leq f(R')$, respectively;
hence $f(A)=f(R')$ and $A=R'$.
\end{proof}

\begin{remark}\label{rmk:pfofovordid}
Whenever the overorders of $R$ are uniquely determined by their conductors,
the last argument gives $R'=\langle R,\pi_E^{f(R')}\Mfo_E\rangle$ for every
overorder $R'$; we use this below.
\end{remark}

\subsubsection{The case of even $e>2$}\label{subsec:3.2.2}

\begin{proposition}\label{prop:charofRin1stcase_1}
Suppose that $e>2$ is even.
\begin{enumerate}
\item There exist $u\in\Mfo_E^\times$ and an odd integer $t\geq3$
such that $R=\Mfo_K[u\pi_E^2,\pi_E^t]$.

\item For odd $k$ with $1\leq k\leq t$, the orders
$\Mfo_k:=\Mfo_K[u\pi_E^2,\pi_E^k]$ are pairwise distinct and exhaust
all overorders of $R$.

\item For each such $k$, $S_K(\Mfo_k)=(k-1)/2$ and
$f(\Mfo_k)=k-1$.
\end{enumerate}
\end{proposition}

\begin{proof}
By Proposition \ref{prop:bassfirstcase}(1), choose
$u\in\Mfo_E^\times$ and put $x:=u\pi_E^2$ and
$A:=\sum_{i=0}^{e/2-1}\Mfo_Kx^i\subset R$. Since $E/K$ is totally
ramified, the elements $x^i$ and $\pi_Ex^i$, $0\leq i<e/2$, have
valuations $2i$ and $2i+1$, respectively, and form an
$\Mfo_K$-basis of $\Mfo_E$. Thus $\Mfo_E=A\oplus A\pi_E$.

Let $t$ be the least valuation of a nonzero element of
$R\cap A\pi_E$, which exists because $R$ is a full $\Mfo_K$-lattice. Every
nonzero element of $A\pi_E$ has odd valuation, so $t$ is odd, and
$u(R)=2$ gives $t\geq3$. Choose $v\in\Mfo_E^\times$ with
$v\pi_E^t\in R\cap A\pi_E$. For $0\leq i<e/2$, the elements
$x^{i+(t-1)/2}$ and $v\pi_E^tx^i$ lie in $R$ and have consecutive
valuations $t-1,\ldots,t+e-2$; hence they form an $\Mfo_K$-basis of
$\pi_E^{t-1}\Mfo_E$. Thus $f(R)\leq t-1$. Since $A\pi_E$ contains
an element of every positive odd valuation, minimality of $t$ gives
$f(R)\geq t-1$. Hence $f(R)=t-1$, and, for $a,b\in A$,
$a+b\pi_E\in R$ if and only if $\ord_E(b\pi_E)\geq t$.

For odd $k$ with $1\leq k\leq t$, the elements
$x^{i+(k-1)/2}$ and $\pi_E^kx^i$, $0\leq i<e/2$, likewise form an
$\Mfo_K$-basis of $\pi_E^{k-1}\Mfo_E$. Hence
$\Mfo_K[x]+\pi_E^{k-1}\Mfo_E\subset\Mfo_k$; the reverse inclusion
holds because the left-hand side is a ring containing $x$ and
$\pi_E^k$. Thus
$\Mfo_k=\Mfo_K[x]+\pi_E^{k-1}\Mfo_E$. For $k=t$, this sum is
contained in $R$. Conversely, write $r=a+b\pi_E\in R$ with
$a,b\in A$. The preceding criterion gives
$b\pi_E\in\pi_E^t\Mfo_E\subset\pi_E^{t-1}\Mfo_E$, while
$a\in A\subset\Mfo_K[x]$. Therefore
$R=\Mfo_K[x]+\pi_E^{t-1}\Mfo_E=\Mfo_t$.

Let $\Mfo$ be an overorder of $R$, and let $k$ be the least valuation
of a nonzero element of $\Mfo\cap A\pi_E$. Then $k$ is odd,
$1\leq k\leq t$, and the preceding argument gives
$\Mfo=\Mfo_k$ and $f(\Mfo_k)=k-1$. Conversely,
$R=\Mfo_t\subseteq\Mfo_k$ for every odd $k$ in this range. Thus the
$\Mfo_k$ are pairwise distinct and exhaust all overorders of $R$;
since $\kappa_{\Mfo_k}=\kappa_E$, Equation
\eqref{eq:relsercon} gives $S_K(\Mfo_k)=(k-1)/2$.
\end{proof}

\begin{theorem}\label{thm:computeorbital39}
Suppose that $e>2$ is even. Then
\[
\left\{
\begin{array}{l}
\#(\Lambda_E\backslash X_R)=1+q^d+\cdots+q^{S(R)};\\
\#\overline{\mathrm{Cl}}(R)=S_K(R)+1
=\frac{S(R)}{d_R}+1=\frac{f(R)}{2}+1.
\end{array}
\right.
\]
Every overorder $R'$ of $R$ satisfies
$R'=\langle R,\pi_E^{f(R')}\Mfo_E\rangle$.
\end{theorem}

\begin{proof}
By Proposition \ref{prop:charofRin1stcase_1}, the overorders are
$\Mfo_k$ for odd $1\leq k\leq t$, with
$S_K(\Mfo_k)=f(\Mfo_k)/2=(k-1)/2$. Since
$\kappa_{\Mfo_k}=\kappa_E$, Proposition \ref{rmk:counting} gives
the weight $q^{d(k-1)/2}$, so Proposition
\ref{prop:charofclbarforbass} gives the first formula by summing and
the second by counting. 
Remark~\ref{rmk:pfofovordid} gives the last assertion.
\end{proof}

\subsection{\texorpdfstring{The case $u(R)=1$ and
$[\kappa_E:\kappa_R]=2$}{The case u(R)=1 and the residue degree is two}}
\label{subsec:2ndcase}
We treat the second case of Proposition \ref{prop:bassclass} and derive
closed formulas for $\#(\Lambda_E\backslash X_R)$ and
$\#\overline{\mathrm{Cl}}(R)$. Assume throughout that
$R\neq\Mfo_E$, $u(R)=1$, and $[\kappa_E:\kappa_R]=2$. Choose a
uniformizer $\pi_E\in R$ with minimal polynomial
$\phi(x)\in\Mfo_K[x]$. Then $[K[\pi_E]:K]$ is either $e$ or $2e$.
In the first case, $K[\pi_E]/K$ is totally ramified of degree $e$ and
$E/K[\pi_E]$ is unramified of degree $2$; in the second,
$K[\pi_E]=E$. We treat the two cases separately below.

\subsubsection{The case $[E:K[\pi_E]]=2$}\label{subsubsec:ek2}
Put $B:=\Mfo_{K[\pi_E]}=\Mfo_K[\pi_E]$. Then $B\subset R$. 

\begin{proposition}\label{prop:charofRin1stcase}
Choose $\alpha\in\Mfo_E^\times$ such that
$\Mfo_E=B\oplus\alpha B$.
\begin{enumerate}
\item There is a unique integer $t\geq1$ such that
$R=B[\alpha\pi_E^t]$.

\item  
The orders $\Mfo_k:=B[\alpha\pi_E^k]$, $0\leq k\leq t$, are pairwise distinct and exhaust all
overorders of $R$.

\item For each such $k$, $S_K(\Mfo_k)=f(\Mfo_k)=k$.
\end{enumerate}
\end{proposition}

\begin{proof}
For an overorder $\Mfo$ of $R$, set
$I_{\Mfo}:=\{a\in B:\alpha a\in\Mfo\}$. Since $B\subset\Mfo\subset\Mfo_E=B\oplus\alpha B$, every element of
$\Mfo$ has a unique expression $b+\alpha a$ with $b\in B$ and
$a\in I_{\Mfo}$; hence $\Mfo=B\oplus\alpha I_{\Mfo}$ and $I_{\Mfo}$ is a nonzero ideal of
the DVR $B$. Thus $I_{\Mfo}=\pi_E^kB$ for a unique
$k\geq0$. Applying this to $R$ gives~(1), with $t\geq1$ because
$R\ne\Mfo_E$. The inclusions $R\subseteq\Mfo\subseteq\Mfo_E$ are
equivalent to $\pi_E^tB\subseteq I_{\Mfo}\subseteq B$, so
$0\leq k\leq t$ and $\Mfo=B\oplus\alpha\pi_E^kB=\Mfo_k$, proving~(2).
The module description gives $S_K(\Mfo_k)=k$. For $k>0$,
$\pi_E^k\Mfo_E\subset\Mfo_k$ but
$\alpha\pi_E^{k-1}\notin\Mfo_k$, so $f(\Mfo_k)=k$; the case $k=0$
is immediate.
\end{proof}

\begin{theorem}\label{thmorb2nd}
Suppose that $[E:K[\pi_E]]=2$. Then
\[
\left\{
\begin{array}{l}
\#(\Lambda_E\backslash X_R)
=q^{S(R)}+2\bigl(q^{S(R)-d_R}+q^{S(R)-2d_R}+\cdots+q^{d_R}+1\bigr);\\
\#\overline{\mathrm{Cl}}(R)=S_K(R)+1
=\frac{S(R)}{d_R}+1=f(R)+1.
\end{array}
\right.
\]
Every overorder $R'$ of $R$ satisfies
$R'=\langle R,\pi_E^{f(R')}\Mfo_E\rangle$.
\end{theorem}

\begin{proof}
By Proposition \ref{prop:charofRin1stcase}, the overorders are
$\Mfo_k$ $(0\leq k\leq t)$, with
$S_K(\Mfo_k)=f(\Mfo_k)=k$. The term for
$\Mfo_0=\Mfo_E$ is $1$, while Proposition \ref{rmk:counting} gives
$(q^{d_R}+1)q^{d_R(k-1)}$ for $1\leq k\leq t$. Proposition
\ref{prop:charofclbarforbass} yields the two formulas by summing and
counting, and Remark~\ref{rmk:pfofovordid} gives the last assertion.
\end{proof}

\subsubsection{The case $E=K[\pi_E]$}\label{subsubsec:ek1}
Assume throughout that $E=K[\pi_E]$. If $e>1$ and $R$ is simple over
$\Mfo_K$, Proposition \ref{prop:polynomialofbass}, applied over $K$,
allows us to choose $\pi_E$ so that
$R=\Mfo_K[\pi_E]$. Let $\widetilde K/K$ be the unramified quadratic
subextension of $E/K$, and choose $u\in\Mfo_{\widetilde K}^\times$ such
that $\Mfo_{\widetilde K}=\Mfo_K\oplus u\Mfo_K$.

For $e>1$, write
\[
\phi(x)=x^{2e}+c_1x^{2e-1}+\cdots+c_{2e-1}x+c_{2e}\in\Mfo_K[x]
\]
for the minimal polynomial of $\pi_E$. Its Newton polygon gives
$\ord_K(c_{2e})=2$, $\ord_K(c_i)\geq1$ for $1\leq i\leq e$, and
$\ord_K(c_i)\geq2$ for $i>e$. Set
$\phi_2(x):=x^2+\overline{c_e/\pi}\,x+\overline{c_{2e}/\pi^2}
\in\kappa_K[x]$.

\begin{lemma}\label{lem:phi_2irrinsep}
The polynomial $\phi_2(x)$ is either irreducible or inseparable.
\end{lemma}

\begin{proof}
Let $\sigma$ be the nontrivial element of
$\operatorname{Gal}(\widetilde K/K)$, and let $g(x)$ be the minimal
polynomial of $\pi_E$ over $\widetilde K$. Since $E/\widetilde K$ is
totally ramified,
$g(x)=x^e+\pi a_1x^{e-1}+\cdots+\pi a_e$ with
$a_e\in\Mfo_{\widetilde K}^\times$, and $\phi=g\,\sigma(g)$. Hence
$\phi_2(x)=(x+\overline{a_e})(x+\sigma(\overline{a_e}))$ over
$\kappa_{\widetilde K}$. If $\overline{a_e}\notin\kappa_K$, the two
roots are distinct conjugates and $\phi_2$ is irreducible; otherwise
they coincide and $\phi_2$ is inseparable.
\end{proof}

Set $\vartheta:=\pi_E^e/\pi$. Reducing
$\phi(\pi_E)/\pi^2$ modulo $\pi_E$ shows that
$\overline\vartheta$ is a root of $\phi_2$. Therefore
\begin{equation}\label{ed:phi2insep}
\overline\vartheta\in\kappa_{\widetilde K}\setminus\kappa_K
\quad\Longleftrightarrow\quad
\phi_2(x)\text{ is irreducible in }\kappa_K[x].
\end{equation}

\begin{proposition}\label{prop:Bassthirdcase}
\begin{enumerate}
\item If $e>1$, $R=\Mfo_K[\pi_E]$, and $\phi_2(x)$ is irreducible,
then $S(R)=d_Re=\frac{n}{2}$.

\item If $e=1$, then $R=\Mfo_K[u\pi^t]$ for a unique integer $t\geq1$.
The orders $\Mfo_k:=\Mfo_K[u\pi^k]$, $0\leq k\leq t$, are pairwise
distinct and exhaust all overorders of $R$, with
$S_K(\Mfo_k)=f(\Mfo_k)=k$.
\end{enumerate}
\end{proposition}

\begin{proof}
For~(1), by \eqref{ed:phi2insep}, the reductions of $1$ and
$\vartheta$ form a $\kappa_K$-basis of $\kappa_E$.
 Reducing the $\Mfo_{\widetilde K}$-basis
$1,\pi_E,\ldots,\pi_E^{e-1}$ of $\Mfo_E$ modulo $\pi$, we see that
$1,\pi_E,\ldots,\pi_E^{e-1},
\pi_E^e/\pi,\ldots,\pi_E^{2e-1}/\pi$ map to a $\kappa_K$-basis of
$\Mfo_E/\pi\Mfo_E$. Nakayama's lemma therefore makes them an
$\Mfo_K$-basis of $\Mfo_E$. Comparison with the power basis of $R=\Mfo_K[\pi_E]$ gives
$S_K(R)=e$, and hence $S(R)=d_Re=n/2$.

For~(2), since $e=1$, $E/K$ is unramified quadratic and
$\Mfo_E=\Mfo_K[u]$. By Remark \ref{rmk:bassideal1}(3), applied over
$K$, every order
$\Mfo$ of $E$ is simple over $\Mfo_K$. Writing a generator as
$a+c\pi^ku$, with $a\in\Mfo_K$, $c\in\Mfo_K^\times$, and $k\geq0$,
gives $\Mfo=\Mfo_K[u\pi^k]=:\Mfo_k$. The basis
$\{1,u\pi^k\}$ gives $S_K(\Mfo_k)=k$, so $k$ is unique. Hence
$R=\Mfo_t$ for a unique $t\geq1$, and its overorders are exactly
 $\Mfo_k$ with $0\leq k\leq t$. Finally,
$\pi^k\Mfo_E\subseteq\Mfo_k$, whereas
$\pi^{k-1}\Mfo_E\nsubseteq\Mfo_k$ for $k>0$; thus
$f(\Mfo_k)=k$, with the case $k=0$ immediate.
\end{proof}

\begin{proposition}\label{prop:insepuniqf}
Suppose that $e>1$, $R=\Mfo_K[\pi_E]$, and $\phi_2(x)$ is inseparable.
Put
$A:=\Mfo_K\oplus\Mfo_K\pi_E\oplus\cdots\oplus
\Mfo_K\pi_E^{e-1}\subset R$ and, for $0\leq k\leq f(R)$, set
$\Mfo_k:=\langle R,\pi_E^k\Mfo_E\rangle=R+\pi_E^k\Mfo_E$.
\begin{enumerate}
\item The orders $\Mfo_k$ are pairwise distinct and exhaust all
overorders of $R$.
\item For each $k$, $S_K(\Mfo_k)=f(\Mfo_k)=k$.
\end{enumerate}
\end{proposition}

\begin{proof}
Since $E/\widetilde K$ is totally ramified of degree $e$,
the powers $1,\pi_E,\ldots,\pi_E^{e-1}$ form an
$\Mfo_{\widetilde K}$-basis of $\Mfo_E$; hence
$\Mfo_E=A\oplus uA$. Since $\phi_2(x)$ is inseparable,
\eqref{ed:phi2insep} gives
$\overline\vartheta\in\kappa_K^\times$.

For $a,b\in A$,
$\ord_E(a+bu)=\min\{\ord_E(a),\ord_E(b)\}$, with
$\ord_E(0)=+\infty$. Indeed, if $\ord_E(a)=\ord_E(b)<+\infty$, write
their common value as $me+r$, $0\leq r<e$, and divide by
$\pi^m\pi_E^r$. The resulting nonzero residues lie in $\kappa_K$ and
$\kappa_K\overline u$, respectively, so they cannot cancel.

Let $\Mfo$ be an overorder of $R$ and set
$k:=\min\{\ord_E(b):a+bu\in\Mfo,\ b\ne0\}$; the set is nonempty
because $\Mfo$ is a full $\Mfo_K$-lattice in $E$. 
Choose $a_0+b_0u\in\Mfo$ with $\ord_E(b_0)=k$. Since
$a_0\in A\subset R$, we have $b_0u\in\Mfo$. Write
$k=me+r$ with $0\leq r<e$. The $2e$ elements
$\pi^m\pi_E^i$ for $r\leq i<e$,
$\pi^{m+1}\pi_E^i$ for $0\leq i<r$, and
$b_0u\pi_E^i$ for $0\leq i<e$ lie in $\Mfo$.

After division by $\pi_E^k$, the first two families together and the
third family each contain one element of every valuation
$j=0,\ldots,e-1$. After dividing the two elements of valuation $j$ by
$\pi_E^j$ and reducing modulo $\pi_E$, their residues lie in
$\kappa_K^\times$ and $\kappa_K^\times\overline u$, respectively.
Here we use $\pi=\vartheta^{-1}\pi_E^e$,
$\overline\vartheta\in\kappa_K^\times$, and
$\overline{b_0/(\pi^m\pi_E^r)}\in\kappa_K^\times$. Since
$\kappa_E=\kappa_K\oplus\kappa_K\overline u$, the two residues at each
valuation are $\kappa_K$-linearly independent. Hence the classes of the $2e$ elements form
a $\kappa_K$-basis of
$\pi_E^k\Mfo_E/\pi\pi_E^k\Mfo_E$.

Nakayama's lemma then gives $\pi_E^k\Mfo_E\subset\Mfo$, so
$f(\Mfo)\leq k$. Conversely, by the definition of $A$, choose
$b\in A$ with $\ord_E(b)=f(\Mfo)$. Since $u$ is a unit,
$bu\in\pi_E^{f(\Mfo)}\Mfo_E\subset\Mfo$, and the minimality of $k$
gives $k\leq f(\Mfo)$. Thus $f(\Mfo)=k$, and
$\Mfo=\{a+bu:a,b\in A,\ \ord_E(b)\geq k\}$.

For $A_j:=A\cap\pi_E^j\Mfo_E$, the identity
$\ord_E(a+bu)=\min\{\ord_E(a),\ord_E(b)\}$ gives
$\pi_E^j\Mfo_E=A_j\oplus uA_j$. The preceding description gives
$R=A\oplus uA_{f(R)}$; since $A_{f(R)}\subseteq A_j$ for
$0\leq j\leq f(R)$, we obtain
$\Mfo_j=R+\pi_E^j\Mfo_E=A\oplus uA_j$.
The minimum defining $k$ for $\Mfo_j$ is $j$, so
$f(\Mfo_j)=j$. Moreover, every overorder $\Mfo$ satisfies
$\Mfo=A\oplus uA_{f(\Mfo)}=\Mfo_{f(\Mfo)}$; hence the $\Mfo_j$ are
pairwise distinct and exhaust all overorders of $R$. For $j>0$,
the above decomposition gives $\kappa_{\Mfo_j}=\kappa_K$, so
\eqref{eq:relsercon} yields $S_K(\Mfo_j)=f(\Mfo_j)=j$; the case
$j=0$ is immediate.
\end{proof}

\begin{theorem}\label{thmorb3rd}
Suppose that $R$ is a simple extension of $\Mfo_K$. If $e>1$, write
$R=\Mfo_K[\pi_E]$ as above and, when $\phi_2(x)$ is irreducible,
assume also that $\operatorname{char}(F)=0$ or
$\operatorname{char}(F)>n$. Then
\[
\#(\Lambda_E\backslash X_R)
=q^{S(R)}+2\bigl(q^{S(R)-d_R}+q^{S(R)-2d_R}
+\cdots+q^{d_R}+1\bigr).
\]
\end{theorem}

\begin{proof}
For $e>1$ with $\phi_2(x)$ irreducible, the formula follows from
Theorem \ref{thmorbital}, Proposition \ref{prop:Bassthirdcase}(1), and
\cite[Lemma~3.2]{CKL}. In the remaining cases, Propositions
\ref{prop:Bassthirdcase}(2) and \ref{prop:insepuniqf} give the
overorders $\Mfo_0,\ldots,\Mfo_s$, where $\Mfo_0=\Mfo_E$,
$S_K(\Mfo_k)=k$, and $s=S_K(R)$. The maximal order has weight $1$,
while Proposition \ref{rmk:counting} gives
$(q^{d_R}+1)q^{d_R(k-1)}$ for $1\leq k\leq s$. The mass identity in
Proposition \ref{prop:charofclbarforbass} then gives the stated formula
because $S(R)=d_Rs$.
\end{proof}

\begin{theorem}\label{thm:overorder2}
Suppose that $R$ is a simple extension of $\Mfo_K$. If $e>1$, write
$R=\Mfo_K[\pi_E]$ as above; if $e=1$, write
$R=\Mfo_K[u\pi^t]$ as in Proposition
\ref{prop:Bassthirdcase}(2). When $e>1$ and $\phi_2(x)$ is
irreducible, assume also that
$\operatorname{char}(F)=0$ or $\operatorname{char}(F)>n$. Set
\[
\Mfo_k:=
\left\{
\begin{array}{ll}
\Mfo_K[\pi_E,\pi_E^k/\pi],
& e\leq k\leq2e,\quad e>1,\ \phi_2(x)\text{ irreducible},\\
\langle R,\pi_E^k\Mfo_E\rangle,
& 0\leq k\leq f(R),\quad e>1,\ \phi_2(x)\text{ inseparable},\\
\Mfo_K[u\pi^k],
& 0\leq k\leq t,\quad e=1.
\end{array}
\right.
\]
\begin{enumerate}
\item The orders $\Mfo_k$ are pairwise distinct and exhaust all
overorders of $R$.
\item
$\#\overline{\mathrm{Cl}}(R)=S_K(R)+1=\frac{S(R)}{d_R}+1
=\left\{
\begin{array}{ll}
e+1,& e>1,\ \phi_2(x)\text{ irreducible},\\
f(R)+1,& e>1,\ \phi_2(x)\text{ inseparable},\\
t+1,& e=1.
\end{array}
\right.$
\end{enumerate}
\end{theorem}

\begin{proof}
For $e=1$ and for $e>1$ with $\phi_2(x)$ inseparable, assertion~(1)
follows from Propositions \ref{prop:Bassthirdcase}(2) and
\ref{prop:insepuniqf}, and assertion~(2) follows by counting and
applying Proposition \ref{prop:charofclbarforbass}. Assume that $e>1$
and $\phi_2(x)$ is irreducible. The basis in Proposition
\ref{prop:Bassthirdcase}(1) shows that the displayed $\Mfo_k$ are
overorders, with $\Mfo_e=\Mfo_E$, $\Mfo_{2e}=R$, and
$S_K(\Mfo_k)=k-e$; hence they are pairwise distinct. The maximal order
has weight $1$, while Proposition \ref{rmk:counting} gives the weight
$(q^{d_R}+1)q^{d_R(k-e-1)}$ for $e+1\leq k\leq2e$. Their sum is the
orbital integral in Theorem \ref{thmorb3rd}, so positivity in the mass
identity proves exhaustion. Counting the candidates and applying
Proposition \ref{prop:charofclbarforbass} give~(2).
\end{proof}

We now remove the simple-extension hypothesis.

\begin{corollary}\label{cor:3rdmain}
Suppose that $E=K[\pi_E]$. When $e>1$ and $\phi_2(x)$ is irreducible,
assume also that $\operatorname{char}(F)=0$ or
$\operatorname{char}(F)>n$. Then
\[
\left\{
\begin{array}{l}
\#(\Lambda_E\backslash X_R)
=q^{S(R)}+2\bigl(q^{S(R)-d_R}+q^{S(R)-2d_R}+\cdots+q^{d_R}+1\bigr);\\
\#\overline{\mathrm{Cl}}(R)=S_K(R)+1=\frac{S(R)}{d_R}+1.
\end{array}
\right.
\]
Every overorder $R'$ of $R$ satisfies
$R'=\langle R,\pi_E^{f(R')}\Mfo_E\rangle$.
\end{corollary}

\begin{proof}
If $e>1$, the uniformizer $\pi_E$ fixed above lies in $R$, so
$R_0:=\Mfo_K[\pi_E]\subseteq R$; it is Bass by
Proposition \ref{prop:bassclass}, since $u(R_0)=1$ and
$[\kappa_E:\kappa_{R_0}]=2$. If $e=1$, take $R_0:=R$, which is simple by Proposition
\ref{prop:Bassthirdcase}(2). Theorem
\ref{thm:overorder2}, applied to $R_0$, shows that the overorders of
$R$ have Serre invariants $0,\ldots,s$, where $s:=S_K(R)$. The
maximal order has weight $1$, and Proposition \ref{rmk:counting} gives
$(q^{d_R}+1)q^{d_R(j-1)}$ for $1\leq j\leq s$. Summing
and counting in Proposition \ref{prop:charofclbarforbass} gives the
two formulas, since $S(R)=d_Rs$.  
Every nonmaximal member of the chain has residue field $\kappa_K$,
while both invariants vanish for $\Mfo_E$. Hence
\eqref{eq:relsercon} gives $f(\Mfo)=S_K(\Mfo)$ throughout the chain;
Remark~\ref{rmk:pfofovordid} gives the final assertion.
\end{proof}

\begin{remark}\label{rmk:anotherpfofsec3}
The local results above use either direct classification or an
exhaustion argument.  In the same-residue-field case
with even $e>2$ (Section~\ref{subsec:3.2.2}),
the residue-degree-two case with inseparable $\phi_2$, and the
unramified quadratic case, the overorders are classified directly.
The orbital integral formulas then follow from the mass identity and
the unit-index formula in Propositions
\ref{prop:orderidealcounting}--\ref{rmk:counting}. 
For odd $e\geq5$ in the same-residue-field case and for irreducible
$\phi_2$ in the residue-degree-two case (see Sections
\ref{sec:first_oddor2} and \ref{subsubsec:ek1}, respectively),  Section \ref{sec:smoothening} computes
the orbital integral used for exhaustion. The same-residue-field cases
$e=2,3$ use the formulas in \cite{CKL}, and the split case uses Yun's
Levi descent formula in \cite[Corollary 4.10]{Yun13}.

Section \ref{sec:smoothening} is therefore a separate geometric result.
It extends the method of \cite{CKL} to the $l_1>0$ strata needed here.
The same construction is globalized in the $\mathrm{GL}_n$ Hitchin
fibration in \cite{CHGlobal}.
\end{remark}

\section{A geometric theorem on smooth integral models and orbital integrals}\label{sec:smoothening}

This section proves Theorem \ref{thmorbital}, which computes the
orbital integral used for exhaustion in the relevant cases of Sections
\ref{sec:first_oddor2} and \ref{subsubsec:ek1}. The proof uses
rescaling and successive dilatations to obtain the required smooth
integral models; following \cite{CKL}, we refer to this procedure as
\emph{smoothening}.
We use the following assumptions throughout this section:
\begin{enumerate}
    \item $\phi(x)$ is an irreducible polynomial of degree $n \geq 4$ in $\mfo[x]$ such that $\ord(\phi(0))=2$.

\item $R\cong \mfo[x]/(\phi(x))$ and $E\cong F[x]/(\phi(x))$ so that $e=\left\{
\begin{array}{l l}
n & \textit{if $n$ is odd};\\
n/2 & \textit{if $n$ is even}.
\end{array}\right.
$ \end{enumerate}
The value of $e$ follows from the Newton polygon of $\phi(x)$. Write
\[
\phi(x)=x^n+c_1x^{n-1}+\cdots+c_{n-1}x+c_n,
\qquad c_i\in\mfo,\quad \ord_F(c_n)=2.
\]
The same Newton polygon gives $\ord_F(c_i)\geq1$ for
$1\leq i\leq\lfloor n/2\rfloor$ and $\ord_F(c_i)\geq2$ for
$\lfloor n/2\rfloor<i\leq n$. When $n$ is even, set
$\phi_2(x):=x^2+\overline{c_{n/2}/\pi}\,x+
\overline{c_n/\pi^2}\in\kappa[x]$. By Lemma
\ref{lem:phi_2irrinsep}, the polynomial $\phi_2(x)$ is either
irreducible or inseparable.
We compute $\#(\Lambda_E\backslash X_R)$ under the following restriction:
\begin{equation}\label{restriction}
\left\{
  \begin{array}{l l}
\phi_2(x) \left(:=x^2+\overline{1/\pi\cdot c_{n/2}}\cdot x+\overline{1/\pi^2\cdot c_{n}}\right) \textit{ is irreducible in }\kappa[x]  & \textit{if $n \geq 4$ is even};\\
\textit{no restriction} & \textit{if $n\geq 5$ is odd}.
    \end{array} \right.
\end{equation}
The cases $n=2,3$ are computed in \cite{CKL}, as recalled in Section \ref{section:orbbass}.
This restriction is used only in Section \ref{subsec:smoothening}, not in Sections \ref{subsec:geom}-\ref{subsec:geomform}.

\begin{remark}\label{rmk:strategywrtgeomeasure}
We now describe the strategy for computing $\#(\Lambda_E\backslash X_R)$.
As explained in Remark \ref{rmk:quotientmeaure}(1),
the number $\#(\Lambda_E \backslash X_R)$ is the same as the orbital integral with respect to the quotient measure.
On the other hand, \cite[Proposition 3.29]{FLN} describes the difference between the quotient measure
and the geometric measure, described precisely in Section \ref{subsec:geom}.
In this context, our method to compute $\#(\Lambda_E \backslash X_R)$ is through the orbital integral with respect to the geometric measure.

We compute the orbital integral for the geometric measure using the stratification in Section \ref{subsec:stratf} and smoothening of certain strata in Sections \ref{subsec:geomform}--\ref{subsec:smoothening}. 
The assumption that $\phi_2(x)$ is irreducible when $n$ is even is used to prove that the relevant stratum is smooth (Proposition \ref{prop:smoothening}).
    Note that this assumption was  previously and independently used in Proposition \ref{prop:Bassthirdcase}(1) to compute $S(R)$.
\end{remark}

\subsection{Orbital integral with respect to the geometric measure}\label{subsec:geom}
The orbital integral of Definition \ref{def:orbitalintegralquo} uses the quotient measure (see Remark \ref{rmk:quotientmeaure}(1)). Here we define the corresponding integral for the geometric measure in Definition \ref{defsoi} and compare the two in Proposition \ref{prop:comp}.

We use the following notation from \cite[Section 2.1]{CKL}.
For a commutative $\mfo$-algebra $A$, let
\[
\left\{
\begin{array}{l}
\textit{$\mathrm{GL}_{n,A}$ be the general linear group of dimension $n^2$ defined over $A$};\\
\textit{$\mathfrak{gl}_{n,A}$ be the Lie algebra of $\mathrm{GL}_{n,A}$};\\
\textit{$\mathbb{A}^{n}_A$ be the affine space of dimension $n$ defined over $A$};\\
\textit{$\chi_{m}(x)\in F[x]$  be the characteristic polynomial of $m \in \mathfrak{gl}_{n, F}(F)$}.
 \end{array}\right.
\]
Let $\omega_{\mathfrak{gl}_{n,\mathfrak o}}$ and $\omega_{\mathbb A^n_{\mathfrak o}}$ be nonzero translation-invariant forms on $\mathfrak{gl}_{n,F}$ and $\mathbb A^n_F$, respectively, normalized by
$\int_{\mathfrak{gl}_{n,\mathfrak o}(\mathfrak o)}
 |\omega_{\mathfrak{gl}_{n,\mathfrak o}}|=1$ and 
$\int_{\mathbb A^n_{\mathfrak o}(\mathfrak o)}
 |\omega_{\mathbb A^n_{\mathfrak o}}|=1$.
 Define $\rho_n:\gl_{n,F}\longrightarrow\mathbb A_F^n$ by
$\rho_n(m)=(c'_1,\ldots,c'_n)$ when
$\chi_m(x)=x^n+c'_1x^{n-1}+\cdots+c'_{n-1}x+c'_n$.
This is a morphism of $F$-schemes.

Note that $\phi(x)$ is viewed as an element of $\mathbb{A}_F^n(F)$.
A differential  $\omega_{\phi}^{\mathrm{ld}}$ on $\rho_n^{-1}(\phi(x))$ is defined to be
$\omega_{\phi}^{\mathrm{ld}}:=\omega_{\mathfrak{gl}_{n, \mathfrak{o}}}/\rho_n^{\ast}\omega_{\mathbb{A}^n_{\mathfrak{o}}}$.
The measure $|\omega_{\phi}^{\mathrm{ld}}|$ is then called \textit{the geometric measure}.
For a precise description of the geometric measure, we refer to \cite[Definition 2.1]{CKL}.

We consider the following morphism of schemes defined over $\mfo$:
$\varphi_n: \mathfrak{gl}_{n, \mathfrak{o}} \longrightarrow \mathbb{A}_{\mathfrak{o}}^n,   ~~~~~~~  m \mapsto
\textit{coefficients of $\chi_{m}(x)$}$.
Thus, the generic fiber of $\varphi_n$ over $F$ is $\rho_n$.
To simplify the notation, let
\begin{equation}\label{desc:ophi}
    O_{\phi}:=\varphi_n^{-1}(\phi)(\mfo)\left(=\{m\in\mathfrak{gl}_{n, \mathfrak{o}}(\mfo)|\chi_m(x)=\phi(x)\}\right).
    \end{equation}

\begin{definition}\label{defsoi}
The  orbital integral for $\phi(x)$  with respect to the geometric measure, denoted by $\mathcal{SO}_{\phi}$, is defined to be
$\mathcal{SO}_{\phi}:=\int_{O_{\phi}}|\omega_{\phi}^{\mathrm{ld}}|$.
\end{definition}

\begin{proposition}(\cite[Propositions 2.4-2.5]{CKL})\label{prop:comp}
Suppose that $\operatorname{char}(F)=0$ or $\operatorname{char}(F)>n$.
Then the following relation holds between $\#(\Lambda_E \backslash X_R)$ and $\mathcal{SO}_{\phi}$:
\[
\#(\Lambda_E \backslash X_R)= \left\{
  \begin{array}{l l}
q^{S(R)}\cdot  \frac{(q-1)q^{n^2-1}}{\#\mathrm{GL}_n(\kappa)}\cdot
\mathcal{SO}_{\phi} & \textit{if $n$ is odd};\\
q^{S(R)}\cdot  \frac{(q^2-1)q^{n^2-2}}{\#\mathrm{GL}_n(\kappa)}\cdot
\mathcal{SO}_{\phi} & \textit{if $n$ is even}.
    \end{array} \right.
\]
\end{proposition}

The following subsections compute $\mathcal{SO}_{\phi}$. Section \ref{subsec:stratf} gives the stratification, Section \ref{subsec:geomform} gives a geometric description of each stratum, and Section \ref{subsec:smoothening} computes the volumes after smoothening.

\subsection{Stratification}\label{subsec:stratf}
We describe four stratifications of $O_\phi$.

\subsubsection{Stratification of $O_\phi$}
Since $\ord_F(\det(m))=\ord_F(c_n)=2$ for $m\in O_{\phi}$, the rank of $\overline{m}$ as a matrix over $\kappa$  is either $n-1$ or $n-2$. If we let
$O^i_\phi=\{m\in O_{\phi}\mid rank(\overline{m})=n-i\}$, then
\begin{equation}\label{equation:ophi12}
    O_{\phi}=O^1_\phi\sqcup O^2_\phi  ~~~ \textit{   so that  }  ~~~~
\mathcal{SO}_{\phi}=\int_{O^1_\phi}|\omega_{\phi}^{\mathrm{ld}}|+\int_{O^2_\phi}|\omega_{\phi}^{\mathrm{ld}}|,
\end{equation}
Note that $O^i_\phi$ is an open subset of $O_{\phi}$ in terms of $\pi$-adic topology.
The volume of $O^1_\phi$ is computed in \cite[Corollary 4.9]{CKL} which is described as follows:
\begin{equation}\label{eq:volo1}
\int_{O^1_\phi}|\omega_{\phi}^{\mathrm{ld}}|=    \frac{\#\mathrm{GL}_{n}(\kappa)}{(q-1)\cdot q^{n^{2}-1}}.
\end{equation}

\subsubsection{Stratification of $O^2_\phi$}
Choose a free $\mfo$-module $L$ of rank $n$ so that we identify  $\mathfrak{gl}_{n, \mfo}(\mfo)=\End_{\mfo}(L)$.
For a sublattice $M$ of $L$, let $O^2_{\phi,M}=\{m\in O^2_\phi\mid m:L\rightarrow M\textit{ is surjective}\}$, which is an open (possibly empty) subset of $O^2_{\phi}$ in terms of $\pi$-adic topology.
If $O^2_{\phi,M}$ is non-empty, then $L/M \cong \mfo/\pi \mfo \oplus \mfo/\pi \mfo$ since $\ord_F(\det(m))=\ord_F(c_n)=2$ and
$rank(\overline{m})=n-2$ for $m\in O^2_{\phi,M}$. Thus we have  the following stratification:
$O^2_{\phi}=\bigsqcup\limits_{L/M\cong \mfo/\pi \mfo \oplus \mfo/\pi \mfo}O^2_{\phi,M}$.
The following formula is proved in \cite[Lemma 3.6 and Corollary 3.9]{CKL}:
   If $L/M \cong L/M'\cong \mfo/\pi \mfo \oplus \mfo/\pi \mfo$, then
\begin{equation}\label{equation:lemma3639}
    \int_{O^2_{\phi,M}}|\omega_{\phi}^{\mathrm{ld}}|=\int_{O^2_{\phi,M'}}|\omega_{\phi}^{\mathrm{ld}}| ~~~~ \textit{ and thus }  ~~~
  \int_{O^2_{\phi}}|\omega_{\phi}^{\mathrm{ld}}|=\frac{(q^n-1)(q^n-q)}{(q^2-1)(q^2-q)}\cdot \int_{O^2_{\phi,M}}|\omega_{\phi}^{\mathrm{ld}}|.
\end{equation}
Here $\frac{(q^n-1)(q^n-q)}{(q^2-1)(q^2-q)}$ is the number of sublattices $M$ of $L$ such that $L/M\cong \mfo/\pi \mfo \oplus \mfo/\pi \mfo$.

\subsubsection{Stratification of $O^2_{\phi,M}$}
We choose a basis $(e_1,\ldots,e_n)$ for $L$ with respect to which
$m\in O^2_{\phi,M}$ is represented by the block matrix
$\begin{pmatrix}
X_{1,1} & X_{1,2} & X_{1,3} \\
\pi X_{2,1} & \pi X_{2,2} & \pi X_{2,3} \\
\pi X_{3,1} & \pi X_{3,2} & \pi X_{3,3}
\end{pmatrix}$,
where $X_{1,1}\in\mathrm{M}_{n-2}(\mfo)$,
$X_{1,j}\in\mathrm{M}_{n-2,1}(\mfo)$ and
$X_{j,1}\in\mathrm{M}_{1,n-2}(\mfo)$ for $j=2,3$, and
$X_{i,j}\in\mfo$ for $i,j=2,3$.
The characteristic polynomial of $\overline{X}_{1,1}$ is
$x^{n-2}\in\kappa[x]$, and its rank is either $n-3$ or $n-4$
since $\ord_F(\det(m))=2$.
Thus there are at most two Jordan blocks for $\overline{X}_{1,1} \left( \in \mathrm{M}_{n-2}(\kappa) \right)$ whose eigenvalues are zero.
We say that $\overline{X}_{1,1}$ is of type $(l_1, l_2)$ if each Jordan block is of size $l_1$ or $l_2$, where $l_1+l_2=n-2$ and $0\leq l_1\leq l_2$.
Then
\begin{equation}\label{eq:strao2fm}
    O^2_{\phi,M}=\bigsqcup\limits^{\lfloor(n-2)/2\rfloor}_{l_1=0}O^2_{\phi, M, l_1}, ~~~~~~  where ~~~~~   O^2_{\phi, M, l_1}=\{m\in O^2_{\phi, M}| ~type(\overline{X}_{1,1})=(l_1, l_2)\}.
\end{equation}

The following formula is proved in \cite[Corollary 7.5]{CKL}:
  When $l_1=0$,  we have
  \begin{equation}\label{equation:l1=0}
         \int_{O^2_{\phi, M, 0}}|\omega_{\phi}^{\mathrm{ld}}|
   =   \frac{\#\mathrm{GL}_{n-2}(\kappa)\cdot (q^{2}-1)}{q^{n^2-2n+3}}.
      \end{equation}

Thus it suffices to analyze the volume of $O^2_{\phi, M, l_1}$ when $l_1>0$.

\subsubsection{Stratification of $O^2_{\phi, M, l_1}$ with $l_1>0$}
Let $J_{l_i}\in\mathrm{M}_{l_i}(\kappa)$ denote the nilpotent Jordan block with superdiagonal entries equal to $1$.
Then $\overline{X}_{1,1}$ is similar to $J_{l_1}\perp J_{l_2}$ if $type(\overline{X}_{1,1})=(l_1, l_2)$.
By the argument used in the proof of \cite[Lemma 3.6]{CKL},
\begin{equation}\label{eq:volo2gln-2}
\int_{O^2_{\phi, M, l_1}}|\omega_{\phi}^{\mathrm{ld}}|=
\#\mathrm{GL}_{n-2}(\kappa)/\#\left(Z_{\mathrm{GL}_{n-2}(\kappa)}(J_{l_1}\perp J_{l_2})\right)\cdot
\int_{\{m\in O^2_{\phi,  M, l_1}|\bar{X}_{1,1}=J_{l_1}\perp J_{l_2}\}}|\omega_{\phi}^{\mathrm{ld}}|.
\end{equation}
Here $\#\mathrm{GL}_{n-2}(\kappa)/\#\left(Z_{\mathrm{GL}_{n-2}(\kappa)}(J_{l_1}\perp J_{l_2})\right)$ is the number of matrices which are conjugate to $J_{l_1}\perp J_{l_2}$.
The denominator is well known by  \cite[Remark (1) in Section 2]{Ful} or \cite[Lemma 2]{Kung}:
\begin{equation}\label{equation:gln-2z}
    \#\left(Z_{\mathrm{GL}_{n-2}(\kappa)}(J_{l_1}\perp J_{l_2})\right)= \left\{
  \begin{array}{l l}
q^{n+2l_1-4}(q-1)^2 & \textit{if $l_2>l_1>0$};\\
q^{2n-7}(q-1)^2(q+1)  & \textit{if $l_2=l_1>0$ (thus $n$ is even)}.
    \end{array} \right.
\end{equation}

\begin{remark}\label{rmk:summarystr}
By stratification explained in  Equations (\ref{equation:ophi12})-(\ref{equation:gln-2z}),
computing $\mathcal{SO}_{\phi}$ reduces to computing
 $\int_{\{m\in O^2_{\phi,  M, l_1}|\bar{X}_{1,1}=J_{l_1}\perp J_{l_2}\}}|\omega_{\phi}^{\mathrm{ld}}|$ with $l_1>0$.
\end{remark}

\subsection{\texorpdfstring{Geometric formulation of $\{m\in O^2_{\phi,M,l_1}\mid\overline{X}_{1,1}=J_{l_1}\perp J_{l_2}\}$ with $l_1>0$}{Geometric formulation for the strata with l1>0}}\label{subsec:geomform}
In this subsection, we express the set $\{m\in O^2_{\phi, M, l_1}|\overline{X}_{1,1}=J_{l_1}\perp J_{l_2}\}$ with $l_1>0$ as the fiber of an algebraic morphism defined over $\mfo$.
Recall from \eqref{desc:ophi} that $O_{\phi}=\varphi_n^{-1}(\phi)(\mfo)$ is the set of matrices $m\in\mathfrak{gl}_{n,\mathfrak{o}}(\mfo)$ with characteristic polynomial $\phi(x)$.
We use the same procedure as in Section \ref{subsec:geom}, with additional congruence conditions on $\mathfrak{gl}_{n,\mathfrak{o}}=\mathrm{End}_\mfo(L)$ and $\mathbb{A}_{\mathfrak{o}}^n$.
In the next subsection, we prove that this fiber is smooth over $\mfo$ and then apply Weil's volume formula \cite[Theorem 2.2.5]{Weil}.

\subsubsection{Geometric formulation  of $\{m\in O^2_{\phi, M, l_1}|\overline{X}_{1,1}=J_{l_1}\perp J_{l_2}\}$}
Define a functor $\mathrm{End}_\mfo(L)_{M,l_1}$ on the category of flat $\mfo$-algebras by
\[
\mathrm{End}_\mfo(L)_{M,l_1}(R)=\left\{m=\begin{pmatrix}
J_{l_1}\perp J_{l_2}+\pi X^\dag_{1,1} &X_{1,2} & X_{1,3} \\
\pi X_{2,1} &\pi X_{2,2} & \pi X_{2,3} \\
\pi X_{3,1} &\pi X_{3,2} & \pi X_{3,3}
\end{pmatrix}\in \mathrm{M}_{n}(R)
\right\},
\]
where the blocks have the sizes specified above and
$m:L\otimes_{\mfo}R\rightarrow M\otimes_{\mfo}R$ is surjective.
This functor is represented by an open subscheme of $\mathbb A^{n^2}_{\mfo}$.
Let $\mathbb A_{l_1}:=\operatorname{Spec}\mfo[r_1,\ldots,r_n]$, and define
\[
s_{l_1}:\mathbb A_{l_1}\longrightarrow\mathbb A^n_{\mfo},\qquad
(r_1,\ldots,r_n)\longmapsto
(\pi r_1,\ldots,\pi r_{l_2+1},
 \pi^2r_{l_2+2},\ldots,\pi^2r_n).
\]
For a flat $\mfo$-algebra $R$, the map $s_{l_1}(R)$ identifies
$\mathbb A_{l_1}(R)=R^n$ with
$(\pi R)^{l_2+1}\times(\pi^2R)^{n-l_2-1}$.
For an arbitrary $\mfo$-algebra $R$, including a $\kappa$-algebra,
the entries of $X^\dag_{1,1}$ and the $X_{i,j}$ are independent
coordinates on $\mathrm{End}_\mfo(L)_{M,l_1}$, while
$\mathbb A_{l_1}(R)=R^n$; the scaled descriptions above apply literally
only when $R$ is flat, as in \cite[p.~802]{C1} and
\cite[p.~467]{C2}. There are unique polynomial functions $r_i$ in the
matrix coordinates such that
\[
\chi_m(x)=x^n+
\sum_{i=1}^{l_2+1}\pi r_i x^{n-i}+
\sum_{i=l_2+2}^{n}\pi^2r_i x^{n-i}.
\]
These functions define an $\mfo$-morphism
$\varphi_{l_1}:\mathrm{End}_\mfo(L)_{M,l_1}\longrightarrow
\mathbb A_{l_1}$.
After base change to $F$, the displayed matrix identifies
$\mathrm{End}_\mfo(L)_{M,l_1}\otimes_\mfo F$ with an open subscheme of
$\mathfrak{gl}_{n,F}$, and $s_{l_1,F}\circ\varphi_{l_1,F}$ is the
restriction of $\rho_n$. The Newton polygon bounds above show that
$(c_1,\ldots,c_n)$ lies in the image of $s_{l_1}(\mfo)$; we again
denote its unique preimage by $\phi$. Hence
\begin{equation}\label{eq:o2fml1-1}
    \{m\in O^2_{\phi, M, l_1}|\overline{X}_{1,1}
    =J_{l_1}\perp J_{l_2}\}
    =\varphi_{l_1}^{-1}(\phi)(\mfo).
\end{equation}

\subsubsection{Comparison between two measures}
Let $\omega_{\mathrm{End}_\mfo(L)_{M,l_1}}$ be the restriction to
$\mathrm{End}_\mfo(L)_{M,l_1}\otimes_\mfo F$ of a nonzero
translation-invariant top form on $\mathbb A^{n^2}_F$ in the matrix
coordinates above, and let $\omega_{\mathbb A_{l_1}}$ be a nonzero
translation-invariant top form on $\mathbb A_{l_1,F}$. Normalize them so
that $\mathbb A^{n^2}_\mfo(\mfo)$ and $\mathbb A_{l_1}(\mfo)$ have
volume $1$, respectively. Following \cite[Definition~2.1]{CKL}, let
$\omega^{ld}_{(\phi,\mathrm{End}_\mfo(L)_{M,l_1})}$ be the quotient of
$\omega_{\mathrm{End}_\mfo(L)_{M,l_1}}$ by
$\varphi_{l_1}^{\ast}\omega_{\mathbb A_{l_1}}$ on
$\varphi_{l_1}^{-1}(\phi)$. Under the generic-fiber identification
above, we have
\[
\begin{gathered}
|\omega_{\mathfrak{gl}_{n,\mathfrak{o}}}|
=|\pi|^{n^2-2n+4}
 |\omega_{\mathrm{End}_\mfo(L)_{M,l_1}}|,\quad
|s_{l_1,F}^{\ast}\omega_{\mathbb A^n_{\mfo}}|
=|\pi|^{2n-l_2-1}|\omega_{\mathbb A_{l_1}}|, \quad
|\omega_{\phi}^{\mathrm{ld}}|
=|\pi|^{n^2-4n+5+l_2}
 |\omega^{ld}_{(\phi,\mathrm{End}_\mfo(L)_{M,l_1})}|.
\end{gathered}
\]
Together with \eqref{eq:o2fml1-1}, this gives:
\begin{proposition}\label{prop52}
\[
\int_{\{m\in O^2_{\phi,  M, l_1}|\bar{X}_{1,1}=J_{l_1}\perp J_{l_2}\}}|\omega_{\phi}^{\mathrm{ld}}|=\int_{\varphi_{l_1}^{-1}(\phi)(\mfo)}|\omega_{\phi}^{\mathrm{ld}}|=q^{-n^2+4n-5-l_2}
\int_{\varphi_{l_1}^{-1}(\phi)(\mfo)}|\omega^{ld}_{(\phi,\mathrm{End}_\mfo(L)_{M, l_1})}|.
\]
\end{proposition}

\subsection{\texorpdfstring{Smoothening of $\varphi_{l_1}^{-1}(\phi)$ with $l_1>0$}{Smoothening for the strata with l1>0}}\label{subsec:smoothening}
By Remark \ref{rmk:summarystr} and Proposition \ref{prop52}, it remains
to compute the volume of $\varphi_{l_1}^{-1}(\phi)(\mfo)$ with respect
to the measure
$|\omega^{ld}_{(\phi,\mathrm{End}_\mfo(L)_{M,l_1})}|$.
Our argument is based on Weil's volume formula stated in
\cite[Theorem 2.2.5]{Weil}.
Namely, if a scheme $S$ is smooth over $\mfo$, then the volume of $S(\mfo)$ with respect to the canonical measure determined by a nowhere-vanishing differential form of top degree is $\#S(\kappa)/q^{\dim(S/\mfo)}$.

Under the restriction (\ref{restriction}) and the hypothesis $\operatorname{char}(F)=0$ or $\operatorname{char}(F)>n$, Proposition \ref{prop:smoothening} proves that $\varphi_{l_1}^{-1}(\phi)$ is smooth over $\mfo$, Corollary \ref{corcounting_1} computes $\#\varphi_{l_1}^{-1}(\phi)(\kappa)$, and Theorem \ref{thmorbital} gives the resulting formula for $\#(\Lambda_E\backslash X_R)$.

\begin{proposition}\label{prop:smoothening}
Suppose that $\operatorname{char}(F)=0$ or $\operatorname{char}(F)>n$. Under the restriction (\ref{restriction}), the scheme $\varphi_{l_1}^{-1}(\phi)$ is smooth over $\mfo$.
\end{proposition}
\proof
By \cite[Theorem~4.6]{CKL}, it is enough to prove that, for every
$m\in\varphi_{l_1}^{-1}(\phi)(\overline{\kappa})$, the differential
$d\varphi_{l_1,m}:T_m\longrightarrow T_{\varphi_{l_1}(m)}$
is surjective. Here $T_m$ and $T_{\varphi_{l_1}(m)}$ are the Zariski
tangent spaces of $\mathrm{End}_\mfo(L)_{M,l_1}\otimes_\mfo
\overline{\kappa}$ and $\mathbb A_{l_1}\otimes_\mfo\overline{\kappa}$
at $m$ and $\varphi_{l_1}(m)$, respectively.
Using the independent matrix coordinates introduced above, write
$m\in\varphi_{l_1}^{-1}(\phi)(\overline{\kappa})$ and $X\in T_m$ as
\begin{equation}\label{eq82}
\left\{
  \begin{array}{l}
m=
\begin{pmatrix}
    J_{l_1} &0&0\\
    0&J_{l_2}&0\\
    0&0&0
\end{pmatrix}+
\begin{pmatrix}
\pi x_{1,1} &\cdots & \pi x_{1,n-2}&x_{1,n-1}&x_{1,n}\\
\vdots  & &\vdots &\vdots &\vdots \\
\pi x_{n-2,1} &\cdots & \pi x_{n-2,n-2}&x_{n-2,n-1}&x_{n-2,n}\\
\pi x_{n-1,1} &\cdots & \pi x_{n-1,n-2}&\pi x_{n-1,n-1}&\pi x_{n-1,n}\\
\pi x_{n,1} &\cdots & \pi x_{n,n-2}&\pi x_{n,n-1}&\pi x_{n,n}
\end{pmatrix};\\
X=\begin{pmatrix}
\pi a_{1,1} &\cdots &\pi a_{1,n-2}&a_{1,n-1}&a_{1,n}\\
\vdots  & &\vdots &\vdots &\vdots \\
\pi a_{n-2,1} &\cdots &\pi a_{n-2,n-2}&a_{n-2,n-1}&a_{n-2,n}\\
\pi a_{n-1,1}  &\cdots & \pi a_{n-1,n-2}&\pi a_{n-1,n-1}&\pi a_{n-1,n}\\
\pi a_{n,1} &\cdots & \pi a_{n,n-2}&\pi a_{n,n-1}&\pi a_{n,n}
\end{pmatrix}.
    \end{array} \right.
\end{equation}
Here $x_{ij},a_{ij}\in\overline{\kappa}$ are the independent
coordinates of the point and the tangent vector on the special fiber,
respectively; the displayed factors of $\pi$ record the scaled matrix
form used above.
For $1\leq i,j\leq n$, let
$E_{ij}\in T_m$ denote the tangent vector with $a_{ij}=1$ and all other
$a_{uv}=0$.
Set
\[
m'=\begin{pmatrix}
    J_{l_1} &0&0\\
    0&J_{l_2}&0\\
    0&0&0
\end{pmatrix}+
\begin{pmatrix}
0 &\cdots &0&x_{1,n-1}&x_{1,n}\\
\vdots &&\vdots&\vdots&\vdots\\
0 &\cdots &0&x_{n-2,n-1}&x_{n-2,n}\\
x_{n-1,1} &\cdots &x_{n-1,n-2}&x_{n-1,n-1}&x_{n-1,n}\\
x_{n,1} &\cdots &x_{n,n-2}&x_{n,n-1}&x_{n,n}
\end{pmatrix}.
\]
Since $\varphi_{l_1}(m)=\phi$, we have
$r_n(m)=\overline{c_n/\pi^2}\in\overline{\kappa}^{\times}$, so $m'$ is
invertible.
Put
\[
Y=\begin{pmatrix}
    x_{l_1,n-1}&x_{l_1,n}\\
    x_{n-2,n-1}&x_{n-2,n}
\end{pmatrix},\qquad
Z=\begin{pmatrix}
    x_{n-1,1}&x_{n-1,l_1+1}\\
    x_{n,1}&x_{n,l_1+1}
\end{pmatrix},\qquad W=YZ.
\]
A block-determinant calculation shows that $m'$ is invertible if and
only if $\det(W)\neq0$.

Let $U\subset T_{\varphi_{l_1}(m)}$ be the subspace on which the
$(l_2+1)$-st and $n$-th coordinates vanish. We first show that $U$ is
contained in the image of $d\varphi_{l_1,m}$, and then treat the two
remaining coordinates.

\begin{enumerate}
\item Suppose that $l_2>l_1$. Since $l_2+1>n/2$, the Newton polygon
bounds give $\ord(c_{l_2+1})\geq2$. 
Thus the $(l_2+1)$-st coordinate of $\varphi_{l_1}(m)$ is zero, and
the characteristic-coefficient calculation identifies it with
$-W_{2,2}$; hence $W_{2,2}=0$. Since $\det(W)\neq0$, we have
$W_{1,2}W_{2,1}\neq0$.

With respect to the coordinates $(dr_1,\ldots,dr_n)$, the images of
$E_{n-2,l}$ for $l_1+1\leq l\leq n-2$ have successive nonzero pivots in
the first $l_2$ coordinates. The images of $E_{l+1,l_1+1}$ for
$0\leq l\leq l_1-1$ have successive pivots $W_{2,1}$ in coordinates
$l_2+2,\ldots,n-1$. These $n-2$ images therefore form a basis of $U$.

Since $\det(W)\neq0$, both $Y$ and $Z$ are invertible. Let
$E_\ast=E_{n-2,n-1}$ if $Z_{1,2}\neq0$, and let
$E_\ast=E_{n-2,n}$ otherwise. Also, let
$E^\ast=E_{l_1,n-1}$ if $Y_{2,2}\neq0$, and let
$E^\ast=E_{l_1,n}$ otherwise. Modulo $U$, the images of $E_\ast$ and
$E^\ast$ in the $(l_2+1)$-st and $n$-th coordinates are, respectively,
$\bigl(-Z_{1,2}\ \text{or }-Z_{2,2},\ast\bigr)
\quad\text{and}\quad
\bigl(0,Y_{2,2}\det(Z)\ \text{or }-Y_{2,1}\det(Z)\bigr)$.
By construction, these pairs are linearly independent. Hence the
induced map to $T_{\varphi_{l_1}(m)}/U$ is surjective, and so is
$d\varphi_{l_1,m}$.

\item Suppose that $l_2=l_1$, so $n=2l_1+2$.
For $P\in\mathrm{GL}_2(\overline{\kappa})$, put
$g_P=\operatorname{diag}(P\otimes I_{l_1},I_2)$. Since $P\otimes I_{l_1}$
commutes with $J_{l_1}\perp J_{l_1}$, conjugation by $g_P$ preserves
$\mathrm{End}_\mfo(L)_{M,l_1}\otimes_\mfo\overline{\kappa}$ and
$\varphi_{l_1,\overline{\kappa}}$, and sends $(Y,Z,W)$ to $(PY,ZP^{-1},PWP^{-1})$. Since $\varphi_{l_1}(m)=\phi$, the
characteristic polynomial of $W$ is $\phi_2$, which is separable by
\eqref{restriction}; hence $W$ is non-scalar. Choosing
$v\in\overline{\kappa}^2$ such that $(v,Wv)$ is a basis and applying
the corresponding conjugation, we may assume $W_{2,1}=1$.

As in (1), the images of $E_{n-2,l}$ for $l_1+1\leq l\leq n-2$ have
successive nonzero pivots in the first $l_1$ coordinates, while those of
$E_{l+1,l_1+1}$ for $0\leq l\leq l_1-1$ have successive pivots
$W_{2,1}=1$ in coordinates $l_1+2,\ldots,n-1$. These $n-2$ images form a
basis of $U$.

Modulo $U$, the remaining coordinates are $-\operatorname{Tr}(W)$ and
$\det(W)$. The images of $E_{l_1,n-1}$ and $E_{l_1,n}$ are
$(-Z_{1,1},Y_{2,2}\det(Z))$ and
$(-Z_{2,1},-Y_{2,1}\det(Z))$, respectively; their determinant is
$\det(Z)\neq0$. Hence the induced map to
$T_{\varphi_{l_1}(m)}/U$ is surjective, and so is
$d\varphi_{l_1,m}$.
\qedhere
\end{enumerate}

\begin{corollary}\label{corcounting_1}
Suppose that $\operatorname{char}(F)=0$ or $\operatorname{char}(F)>n$. Under the restriction (\ref{restriction}), the cardinality of the set  $\varphi_{l_1}^{-1}(\phi)(\kappa)$ is
\[
\#\varphi_{l_1}^{-1}(\phi)(\kappa)=(q-1)q^{n^2-n-5}\cdot \#\mathrm{GL}_2(\kappa).
\]
\end{corollary}

\begin{proof}
Write $\phi(x)$, viewed as a point of $\mathbb A_{l_1}(\mfo)$, as
\begin{equation}\label{eq:fx}
\phi(x)=x^n+\pi k_1x^{n-1}+\cdots+\pi k_{l_2+1}x^{n-l_2-1}
+\pi^2k_{l_2+2}x^{n-l_2-2}+\cdots+\pi^2k_{n-1}x+\pi^2k_n,
\end{equation}
where $k_i\in\mfo$, $k_n\in\mfo^\times$, and $\overline{k_i}$ denotes
the image of $k_i$ in $\kappa$. If $l_2>l_1$, then
$l_2+1>n/2$, so the Newton polygon bounds give
$k_{l_2+1}\in\pi\mfo$.
For $m\in\varphi_{l_1}^{-1}(\phi)(\kappa)$, write the coordinate matrix
$(x_{ij})$ in \eqref{eq82} as
$(x_{ij})=
\begin{pmatrix}
X_1&X_2\\
X_3&X_4
\end{pmatrix}$,
where $X_1\in\mathrm M_{n-2}(\kappa)$,
$X_2\in\mathrm M_{n-2,2}(\kappa)$,
$X_3\in\mathrm M_{2,n-2}(\kappa)$, and
$X_4\in\mathrm M_2(\kappa)$. Then $Y$ and $Z$ are the indicated
$2\times2$ submatrices of $X_2$ and $X_3$.

\begin{itemize}
    \item
The equations for the $(l_2+1)$-st and $n$-th coordinates of $\varphi_{l_1}(m)$ are
\begin{equation}\label{eq:a1}
    \left\{
  \begin{array}{l l}
W_{2,2}=0, ~~~ \det(W)=\overline{k_n} & \textit{if $l_2>l_1$};\\
-\mathrm{Tr}(W)=\overline{k_{l_2+1}}, ~~~ \det(W)=\overline{k_n} & \textit{if $l_2=l_1$ (thus $n$ is even)}.
    \end{array} \right.
\end{equation}
These equations involve only the entries of $Y$ and $Z$.

\item For $1\leq i\leq l_2$, the $i$-th coordinate of $\varphi_{l_1}(m)$ gives
\begin{equation}\label{eq:a2}
-x_{n-2, n-1-i} + \ast=\overline{k_i}, \qquad\textit{where $\ast$ does not involve $x_{n-2, n-1-i}$.}
\end{equation}

\item For $l_2+2\leq i\leq n-1$, the $i$-th coordinate of $\varphi_{l_1}(m)$ gives
\begin{equation}\label{eq:a3}
-x_{i-l_2-1, l_1+1}\cdot W_{2,1}+\ast=\overline{k_i}, \qquad \textit{where $\ast$ does not involve $x_{i-l_2-1,l_1+1}$.}
\end{equation}
\end{itemize}

The variables singled out in \eqref{eq:a2} occur neither in
\eqref{eq:a1} nor in the other equations in \eqref{eq:a2}. For
$i<j$, the variable singled out in the $i$-th equation of
\eqref{eq:a3} occurs neither in the $j$-th equation nor in
\eqref{eq:a1}--\eqref{eq:a2}. Hence \eqref{eq:a2} determines $l_2$
coordinates and, if $W_{2,1}\neq0$, the equations in \eqref{eq:a3}
for $i=n-1,n-2,\ldots,l_2+2$ determine $l_1$ more, for a total of
$n-2$.

If $l_2>l_1$, \eqref{eq:a1} gives $W_{2,2}=0$ and
$-W_{1,2}W_{2,1}=\overline{k_n}\neq0$. Hence there are $q(q-1)$
choices for $W$, with $W_{1,1}\in\kappa$ and
$W_{1,2}\in\kappa^\times$ free. If $l_2=l_1$, then $W$ has
irreducible characteristic polynomial $\phi_2$, so $W_{2,1}\neq0$.
The possible $W$ form a single $\mathrm{GL}_2(\kappa)$-conjugacy class
whose centralizer in $\mathrm{GL}_2(\kappa)$ has order $q^2-1$; hence
there are again $q(q-1)$ choices.

For each $W$, choosing $Y\in\mathrm{GL}_2(\kappa)$ determines
$Z=Y^{-1}W$. The remaining
$n^2-8-(n-2)=n^2-n-6$ coordinates are free. Hence the total number is
$(q-1)q^{n^2-n-5}\#\mathrm{GL}_2(\kappa)$, as claimed. 
\end{proof}

\begin{remark}\label{rmk:reasonofrestriction}
In Proposition \ref{prop:smoothening} and Corollary \ref{corcounting_1}, the case $l_2>l_1$ does not require \eqref{restriction}. The restriction is used only when $l_2=l_1$, which can occur only for even $n$. Thus the restriction (\ref{restriction}) imposes no additional condition when $n$ is odd.
\end{remark}

\begin{theorem}\label{thmorbital}
Suppose that $\operatorname{char}(F)=0$ or $\operatorname{char}(F)>n$.
For $n=2m+1$ or $n=2m$, the following formula holds:
\[
\#(\Lambda_E \backslash X_R)=
\left\{
   \begin{array}{l l}
q^{S(R)}+q^{S(R)-1}+\cdots +q^{S(R)-m} & \textit{if $n=2m+1\geq 5$};\\
q^{S(R)}+2\left(q^{S(R)-1}+\cdots +q^{S(R)-m}\right)  & \textit{if  $n=2m\geq 4$ and $\phi_2$ irreducible}.
     \end{array} \right.
\]

\end{theorem}
\begin{proof}
The preceding smoothness and point-counting results, together with Weil's
volume formula \cite[Theorem~2.2.5]{Weil}, give
$\int_{\varphi_{l_1}^{-1}(\phi)(\mfo)}
|\omega^{ld}_{(\phi,\mathrm{End}_\mfo(L)_{M,l_1})}|
=\frac{\#\varphi_{l_1}^{-1}(\phi)(\kappa)}{q^{n^2-n}}
=(q-1)q^{-5}\#\mathrm{GL}_2(\kappa)$.
Proposition \ref{prop52} then gives
$\int_{\varphi_{l_1}^{-1}(\phi)(\mfo)}
|\omega_{\phi}^{\mathrm{ld}}|
=q^{-n^2+3n-8+l_1}(q-1)\#\mathrm{GL}_2(\kappa)$.
Using Proposition \ref{prop:comp} and
Equations \eqref{equation:ophi12}--\eqref{equation:gln-2z}
(see Remark \ref{rmk:summarystr}), we obtain:
\begin{enumerate}
    \item Suppose that $n=2m+1$. Then
\begin{multline*}
\#(\Lambda_E \backslash X_R)=q^{S(R)}
+q^{S(R)}\cdot  \frac{q^{n^2-2}(q^n-1)(q^n-q)}{(q^2-1)\cdot\#\mathrm{GL}_n(\kappa)}\cdot \\
  \left(    \frac{\#\mathrm{GL}_{n-2}(\kappa)\cdot (q^{2}-1)}{q^{n^2-2n+3}}+\sum_{l_1=1}^{m-1}
  \frac{\#\mathrm{GL}_{n-2}(\kappa)}{(q-1)^2q^{n+2l_1-4}}\cdot q^{-n^2+3n-8+l_1}\cdot (q-1)\cdot \#\mathrm{GL}_2(\kappa) \right)\\
=  q^{S(R)} +q^{S(R)-1}+
 q^{S(R)}\cdot  \frac{(q^n-1)(q^n-q)}{(q^2-1)\cdot\#\mathrm{GL}_n(\kappa)}\cdot \left(   \sum_{l_1=1}^{m-1}\frac{\#\mathrm{GL}_{n-2}(\kappa)}{(q-1)q^{n+2l_1-4}}\cdot q^{3n-10+l_1}\cdot \#\mathrm{GL}_2(\kappa)
  \right).
  \end{multline*}
To summarize,
$\#(\Lambda_E \backslash X_R)=
  q^{S(R)} +q^{S(R)-1}+
 q^{S(R)}\cdot \frac{q^{m-1}-1}{q^m(q-1)}=q^{S(R)} +q^{S(R)-1}+\cdots + q^{S(R)-m}$.
Here, we use the formula
$\#\mathrm{GL}_n(\kappa)=\#\mathrm{GL}_{n-2}(\kappa)\cdot (q^n-1)  (q^{n}-q)q^{2n-4}.$

   \item Suppose that $n=2m$ is even and that $\phi_2(x)$ is irreducible in $\kappa[x]$.
 Then
\begin{multline*}
\#(\Lambda_E \backslash X_R)=(q+1)q^{S(R)-1}+(q+1)q^{S(R)-2}
+q^{S(R)-1}\cdot  \frac{q^{n^2-2}(q^n-1)(q^n-q)(q+1)}{(q^2-1)\cdot\#\mathrm{GL}_n(\kappa)}\cdot \\
  \left(
\frac{\#\mathrm{GL}_{n-2}(\kappa)\cdot \#\mathrm{GL}_{2}(\kappa)\cdot (q-1)q^{-n^2+3n-9+m}}
  {(q-1)^2(q+1)q^{2n-7}}
  \quad+\sum_{l_1=1}^{m-2}\frac{\#\mathrm{GL}_{n-2}(\kappa)}{(q-1)^2q^{n+2l_1-4}}
  q^{-n^2+3n-8+l_1}(q-1)\#\mathrm{GL}_2(\kappa)
  \right)\\
=  (q+1)q^{S(R)-1}+(q+1)q^{S(R)-2}
+q^{S(R)-1}(q+1)\cdot
  \left( \frac{q^{1-m}}{q+ 1}  +\sum_{l_1=1}^{m-2}  q^{-1-l_1}
  \right).
\end{multline*}
\end{enumerate}

To summarize,
$\#(\Lambda_E \backslash X_R)
=q^{S(R)}+2\left(q^{S(R)-1}+\cdots +q^{S(R)-m}\right).$
\end{proof}

\part{Global Theory}\label{part2}

Part 2 gives a conductor formula for $\#\left(\mathrm{Cl}(R)\backslash\overline{\mathrm{Cl}}(R)\right)$ when $R$ is a Bass order in a number field.

\section*{Notations}\label{section:globalnotations}

We reset the following notations.
\begin{itemize}
\item For a ring $A$, the set of maximal ideals is denoted by $|A|$.
\item For a local ring $A$, the maximal ideal is denoted by $\mathfrak{m}_A$, and the residue field is denoted by $\kappa_A$. If $\kappa_A$ is finite, put $q_A=\#\kappa_A$.
If $K$ is  a non-Archimedean local field, then  we sometimes use $\kappa_K$ to denote the residue field of the ring of integers in $K$, if there is no confusion.

\item If $A$ is a local ring, then $\overline{a}\in\kappa_A$ and
$\overline{\psi(x)}\in\kappa_A[x]$ denote the reductions of
$a\in A$ and $\psi(x)\in A[x]$ modulo $\mathfrak m_A$, respectively.

\item Let $F$ be a number field with $\mfo$ its ring of integers.
\item For $v\in |\mfo|$,
let $F_v$ be the $v$-adic completion of $F$ with $\mfo_v$  the ring of integers of $F_v$, $\pi_v$  a uniformizer in $\mfo_v$, and $\kappa_v$ its residue field. Let $q_v=\#\kappa_v$.

\item Let $E/F$ be a finite field extension and $R$ an order in $E$ as in
Definition A.(1) below.
\end{itemize}

Following \cite{Ma24}, we recall the relevant notions for orders in
finite \'etale algebras. This setting is used in the local--global
argument and in the split case of Section
\ref{subsection:formulaforRsplit}.

\begin{defff}{(cf. \cite[Sections 2.1--2.3]{Ma24})}
Let $Z$ be a Dedekind domain with field of fractions $Q$.
Let $K$ be an \'etale $Q$-algebra.
Note that $Z$ is always  $\mfo$ or $\mfo_v$ in this paper.

\begin{enumerate}
\item{\cite[the second paragraph of Section 2.2]{Ma24}}  An order of $K$ is a subring $\Mfo$ of $K$ such that $\Mfo$ is a finitely generated $Z$-module containing $Z$ and such that $\Mfo\otimes_ZQ\cong K$.

\item{\cite[the first paragraph of Section~2.3]{Ma24}} A fractional
$\Mfo$-ideal is a finitely generated $\Mfo$-submodule $I$ of $K$ with
$I\otimes_ZQ\cong K$. The same reference gives the following criterion.

\begin{enumerate}
    \item If $I$ is a finitely generated $\Mfo$-submodule of $K$ (e.g.  an ideal of $\Mfo$), then $I$ is a fractional $\Mfo$-ideal if and only if $I$ contains a non-zero divisor of $K$.

\item If $\Mfo$ is a Noetherian domain, then the above criterion yields that a finitely generated  $\Mfo$-submodule $I (\neq 0)$ of $K$ is a fractional $\Mfo$-ideal.
In particular, if $xI\subset \Mfo$ for an $\Mfo$-submodule $I (\neq 0)$ of $K$ with   $x\neq 0$ in  $\Mfo$, then $I$ is a fractional $\Mfo$-ideal.
\end{enumerate}
\item Let $\Mfo_K$ denote the unique maximal order of $K$; see
\cite[the first paragraph of Section~2]{Ma20} or
\cite[the third paragraph of Section~2.2]{Ma24}.

\item(Generalization of Definition \ref{def:invariantofO}(2)) The conductor $\mathfrak{f}(\Mfo)$ of an order $\Mfo$ is the largest ideal of $\Mfo_K$ which is contained in $\Mfo$. In other words, $\mathfrak{f}(\Mfo)=\{a\in \Mfo_K\mid a\Mfo_K\subset \Mfo\}$.
Note that $\mathfrak{f}(\Mfo)$ is also an ideal of $\Mfo$.

\item We use $\langle\Mfo,I\rangle$, ideal quotients, invertible
$\Mfo$-ideals, $\mathrm{Cl}(\Mfo)$, and
$\overline{\mathrm{Cl}}(\Mfo)$ as in the Notations of Part
\ref{part1}. For $[I]\in\overline{\mathrm{Cl}}(\Mfo)$, the
multiplicator ring $(I:I)$ is the largest order over which $I$ is a
fractional ideal. We refer to
\cite[Sections~2.1--2.3]{Ma24} for details.
\end{enumerate}
\end{defff}
\[
\textit{For $v\in |\mfo|$, let  }
\left\{
\begin{array}{l}
\textit{$R_v\cong
 R\otimes_\mfo \mfo_v$ be the $v$-adic completion of $R$};\\
\textit{$E_v\cong R_v\otimes_{\mfo_v} F_v$ be the ring of total fractions of $R_v$};\\
\textit{$X_{R_v}$ be the set of fractional $R_v$-ideals so that $\overline{\cl}(R_v)=E_v^\times\backslash X_{R_v}$}.
\end{array} \right.
\]
\[
\textit{    For $w\in |R|$, let }
\left\{
\begin{array}{l}
\textit{$R_w$ be the $w$-adic completion of $R$};\\
\textit{$E_w\cong E\otimes_R R_w$ be the ring of total fractions of $R_w$};\\
\textit{$X_{R_w}$ be the set of fractional $R_w$-ideals  so that $\overline{\cl}(R_w)=E_w^\times\backslash X_{R_w}$}.
\end{array} \right.
\]
Each $R_w$ is local, but need not be a domain. For $w\mid v$, let
$K_w\subset E_w$ be the unramified field extension of $F_v$ corresponding to
the residue field extension $\kappa_{R_w}/\kappa_v$. The proof of
\cite[Lemma~3.1]{CKL} shows that $\Mfo_{K_w}\subset R_w$; we use $K_w$
in Section \ref{sec:formulabass}.

\begin{remmm}\label{rmk:descriptionofRv}
This remark records the relations among $R_v$, $R_w$, $E_v$, $E_w$, $\Mfo_{E_v}$, and $\Mfo_{E_w}$.
    In (2)--(4), we work with $R\cong \mfo[x]/(\phi(x))$ for an irreducible polynomial $\phi(x)\in \mfo[x]$.

    \begin{enumerate}
\item By \cite[(4) in the proof of Lemma 2.16]{Ma24}, we have
$R_v\cong \bigoplus\limits_{w|v,~ w\in |R|}R_w$.
Applying $(-)\otimes_{\mfo_v}F_v$ gives
$E_v\cong\bigoplus_{w\mid v}E_w$, and hence
$\Mfo_{E_v}\cong\bigoplus_{w\mid v}\Mfo_{E_w}$.
 Note that $E_w$ may not be a field.

\item Consider $\phi(x)$ as an element of $\mfo_v[x]$ for $v\in |\mfo|$ and write
    $\overline{\phi(x)}=(\overline{g_1(x)})^{n_1}\cdots (\overline{g_r(x)})^{n_r}$, where $\overline{g_i(x)}$'s are distinct and irreducible in $\kappa_v[x]$.
    Then    Hensel's lemma yields a factorization $\phi(x)=\phi_1(x)\cdots \phi_r(x)$ in $\mfo_v[x]$ such that $\overline{\phi_i(x)}=(\overline{g_i(x)})^{n_i}$.
We have
        \begin{equation}\tag{A1}\label{eq:decomp_compl}
R_v\cong \prod\limits_{i=1}^r \mfo_v[x]/(\phi_i(x))\cong \bigoplus_{w|v,~ w\in |R|}R_w,
\end{equation}
where  $R_w\cong \mfo_v[x]/(\phi_i(x))$ for some $i$, which is a local ring.
Here the first isomorphism follows from the Chinese remainder theorem, and the second follows from (1).

\item Let $\phi_w(x)$ denote the factor $\phi_i(x)$ corresponding to $R_w$; thus
$R_w\cong\mfo_v[x]/(\phi_w(x))$.
The polynomial $\phi_w(x)$ need not be irreducible, so $R_w$ need not be an integral domain.  Nevertheless, $R_w$ is reduced because $\phi(x)$ is separable over $F_v$; see \cite[Lemma 2.16]{Ma24}.

        \item
Let $B(\phi_w)$ be the index set of irreducible factors of
$\phi_w(x)$ in $\mfo_v[x]$. Then
        \[
        E_v\cong \bigoplus\limits_{w|v, \ w\in|R|} F_v[x]/(\phi_w(x))\cong\bigoplus\limits_{w|v,~ w\in |R|}E_w ~~~~  \textit{  and   } ~~~~~  E_w\cong F_v[x]/(\phi_w(x)) \cong\bigoplus\limits_{j\in B(\phi_w)} E_{w,j},
        \]
        where $E_{w,j}$ is a finite field extension of  $F_v$.
        Here the first follows from  Equation (\ref{eq:decomp_compl}). Thus
\[
\Mfo_{E_v}\cong \bigoplus\limits_{w\mid v, ~ w\in |R|} \Mfo_{E_w} ~~~~~~~  \textit{ and  } ~~~~~~~~~~~~ \Mfo_{E_w}\cong \bigoplus\limits_{j\in B(\phi_w)} \Mfo_{E_{w,j}}.
\]
Note that the Chinese Remainder Theorem implies the inclusion $R_w \hookrightarrow \bigoplus\limits_{j\in B(\phi_w)} \Mfo_{E_{w,j}}$.

    \end{enumerate}
\end{remmm}

\section{Ideal class monoids: local-global argument}\label{sec:uppbddsimple}
This section first stratifies $\overline{\mathrm{Cl}}(\Mfo)$ and $\mathrm{Cl}(\Mfo)\backslash\overline{\mathrm{Cl}}(\Mfo)$ for an order $\Mfo$ in the setting of Definition A, and then proves the local--global statement for $\mathrm{Cl}(R)\backslash\overline{\mathrm{Cl}}(R)$.

\begin{definition}[Generalization of Definition \ref{def:clmfo}]\label{def:globalclcl}
Let $\mathcal{O}$ be an order in the sense of Definition A.
For an overorder $\mathcal{O}'$ of $\mathcal{O}$, we define the following sets:
\[\left\{
\begin{array}{l}
     \mathrm{cl}(\mathcal{O}'):=\{ [I] \in \overline{\mathrm{Cl}}(\mathcal{O}') \mid (I:I)=\mathcal{O}' \}=\{ [I] \in \overline{\mathrm{Cl}}(\mathcal{O}) \mid (I:I)=\mathcal{O}' \};  \\
     \overline{\mathrm{cl}(\mathcal{O}')}:=\mathrm{Cl}(\Mfo)\backslash \mathrm{cl}(\mathcal{O}').
\end{array}\right.
\]
Here we consider $\overline{\mathrm{Cl}}(\mathcal{O}')$ as a subset of $\overline{\mathrm{Cl}}(\mathcal{O})$.
The set  $\overline{\mathrm{cl}(\mathcal{O}')}$ is well-defined since $(JI:JI)=JJ^{-1}(I:I)=(I:I)$ for $J\in \mathrm{Cl}(\mathcal{O})$.
Note that $\mathrm{cl}(\mathcal{O}')$ is denoted by $\mathrm{ICM}_{\Mfo'}(\Mfo)$ in \cite{Ma24}.

\end{definition}

\begin{proposition}\label{prop:stratforglobal}
For an order $\Mfo$ in the sense of Definition A,
the following statements hold:
\[
\#\overline{\mathrm{Cl}}(\mathcal{O})=\sum_{\mathcal{O}\subset \mathcal{O}'\subset \mathcal{O}_K}\#\mathrm{cl}(\mathcal{O}')
~~~~~~~~~\textit{ and } ~~~~~~~~~~~
\mathrm{Cl}(\mathcal{O})\backslash \overline{\mathrm{Cl}}(\mathcal{O})
=\bigsqcup_{\mathcal{O}\subset \mathcal{O}'\subset \mathcal{O}_K}\overline{\mathrm{cl}(\mathcal{O}')}.
\]
\end{proposition}
\proof
These two are direct consequences of the following stratification:
    \[
    \overline{\mathrm{Cl}}(\Mfo)=\bigsqcup_{\Mfo\subset \mathcal{O}' \subset \mathcal{O}_{K}}\{ [I]\in \overline{\mathrm{Cl}}(\Mfo)\mid (I:I)=\mathcal{O}' \}=\bigsqcup_{\Mfo\subset \mathcal{O}' \subset \mathcal{O}_{K}}\mathrm{cl}(\mathcal{O}'). \qedhere
    \]

The following proposition proves the local--global statement for $\mathrm{Cl}(R)\backslash \overline{\mathrm{Cl}}(R)$.
Note that this is mentioned in \cite[Section 3.4, p. 410]{Yun13} when $R$ is a simple extension of $\mfo$, without a proof.

\begin{proposition}\label{prop:global_local}
For an order $R$ of a number field $E$, the following map is bijective:
\[
\mathrm{Cl}(R)\backslash \overline{\mathrm{Cl}}(R)\longrightarrow
\prod_{w\in |R|}\overline{\cl}(R_w),
\qquad
\{I\}\longmapsto \bigl([I\otimes_RR_w]\bigr)_{w\in|R|}.
\]
\end{proposition}

\begin{proof}
By \cite[Lemma~2.17]{Ma24}, the extension of every invertible
fractional $R$-ideal to each $R_w$ is principal; hence the map is
well-defined.

To prove surjectivity, let $([I_w])_{w\in|R|}$ be an element of the
target and choose representatives $I_w$, with $I_w=R_w$ for almost
all $w$. For $v\in|\mfo|$, set $I_v:=\prod_{w\mid v}I_w$ via the
decompositions in Remark B.(1). Then $I_v$ is a fractional
$R_v$-ideal and $I_v=R_v$ for almost all $v$. Since $I_v$ and $R_v$
are full lattices in $E_v$ over the DVR $\mfo_v$, choose
$n_v\geq0$, with $n_v=0$ for almost all $v$, such that
$\pi_v^{n_v}R_v\subset I_v\subset\pi_v^{-n_v}R_v$. Choose
$0\neq a\in\prod_{n_v>0}v^{n_v}$; then
$aR_v\subset I_v\subset a^{-1}R_v$ for every $v$.

Let $S=\{v\in|\mfo|:v\mid a\mfo\}$, and write
$a\mfo=\prod_{v\in S}v^{e_v}$. Since
$a^2R=\prod_{v\in S}v^{2e_v}R$ with pairwise comaximal factors, and
$R/v^{2e_v}R\cong R_v/a^2R_v$ for $v\in S$, the Chinese remainder
theorem gives
\[
R/a^2R\cong\prod_{v\in S}R/v^{2e_v}R
\cong\prod_{v\in S}R_v/a^2R_v.
\]
Multiplying by $a^{-1}$ gives the natural $R$-linear isomorphism
\[
\Phi:a^{-1}R/aR\xrightarrow{\sim}
\prod_{v\in S}a^{-1}R_v/aR_v,
\qquad
x+aR\longmapsto (x+aR_v)_{v\in S}.
\]
Let $I\subset a^{-1}R$ be the inverse image under $\Phi$ of
$\prod_{v\in S}I_v/aR_v$. Then $aR\subset I\subset a^{-1}R$, so $I$
is a full $\mfo$-lattice. Since $\Phi$ is $R$-linear and each $I_v$
is an $R_v$-ideal, $I$ is $R$-stable and hence a fractional
$R$-ideal. For $v\in S$, flatness of $\mfo_v$ over $\mfo$ preserves
this inverse-image construction, while tensoring with $\mfo_v$ kills
the factors indexed by $u\ne v$. Thus
$I\otimes_\mfo\mfo_v\cong I_v$. If $v\notin S$, then $a$ is a unit
in $\mfo_v$, so both sides are equal to $R_v$. Taking the
$w$-component gives $I\otimes_RR_w\cong I_w$ for every $w$, proving
surjectivity.

To prove injectivity, suppose that $I$ and $J$ have the same image.
Choose $\alpha_w\in E_w^\times$ with
$J\otimes_RR_w=\alpha_w(I\otimes_RR_w)$ for every $w$ and
$\alpha_w=1$ for almost all $w$, and put
$L_v=\prod_{w\mid v}\alpha_wR_w$. The preceding construction gives a
fractional $R$-ideal $L$ with $L\otimes_RR_w=\alpha_wR_w$ for every
$w$; it is invertible by \cite[Lemma~2.17]{Ma24}. Since
$(LI)\otimes_RR_w=J\otimes_RR_w$ for every $w$, the
$\mfo_v$-completions of $LI$ and $J$ agree for every $v$; hence
$LI=J$, as full $\mfo$-lattices in $E$ are determined by their
completions.
\end{proof}

The following corollary shows that an overorder of $R$ is determined by its completions at all $w\in|R|$; this fact is crucial in deriving the formula for $\#\left(\mathrm{Cl}(R)\backslash\overline{\mathrm{Cl}}(R)\right)$ in the next section.

\begin{corollary}\label{cor:prop_global_local}
\begin{enumerate}
    \item
    The following mapping is a bijection,
    \[
    \{\textit{overorders $\mathcal{O}$ of $R$}\} \cong \prod_{w\in|R|} \{\textit{overorders $\mathcal{O}_w
$ of $R_w$}\},\ \mathcal{O}\mapsto \prod\limits_{w\in|R|}(\mathcal{O}\otimes_{R}R_w).
    \]
    \item
    For an overorder $\mathcal{O}$ of $R$, we have the bijection,
    \[\overline{\mathrm{cl}(\mathcal{O})}\cong\prod\limits_{w\in|R|} \mathrm{cl}(\mathcal{O}\otimes_R R_w), \ \{I\}\mapsto \prod_{w\in |R|} [I\otimes_{R}R_w].\]
    \item
    The action of $\mathrm{Cl}(\mathcal{O})$ on $\mathrm{cl}(\mathcal{O})$  is free.
    In addition, we have the following identities:
    \[
    \#\mathrm{cl}(\mathcal{O})=\#\mathrm{Cl}(\mathcal{O})\cdot\#\overline{\mathrm{cl}(\mathcal{O})}=\#\mathrm{Cl}(\mathcal{O})\cdot \prod_{w\in|R|}
    \#\mathrm{cl}(\mathcal{O}\otimes_R R_w).
    \]
\end{enumerate}
\end{corollary}
\begin{proof}
For a fractional $R$-ideal $I$, 
put $I_w=I\otimes_RR_w$. Since $I$ is finitely presented and $R_w$
is flat over $R$, \cite[Lemma~2.1 and Equation~(2)]{Ma24} and 
\cite[\href{https://stacks.math.columbia.edu/tag/087R}{Lemma 087R}]
{stacks-project} give
$(I:I)\otimes_RR_w\cong(I_w:I_w)$.
\begin{enumerate}
\item
The map is well-defined. Given $(\mathcal O_w)_w$ in the target,
Proposition \ref{prop:global_local} gives an orbit $\{I\}$ with
$[I_w]=[\mathcal O_w]$ for every $w$. Hence
$\mathcal O=(I:I)$ has completions $\mathcal O_w$, proving
surjectivity. If two overorders have the same completions, Proposition
\ref{prop:global_local} places them in the same $\mathrm{Cl}(R)$-orbit;
they are equal because multiplicator rings are constant on such
orbits.

\item This follows from Proposition \ref{prop:global_local}, the preceding base-change identification, and~(1).

\item
The action is free by \cite[Theorem~4.6]{Ma20}; together with~(2), this gives the identities. \qedhere
\end{enumerate}
\end{proof}

\section{\texorpdfstring{A formula for $\#\left(\mathrm{Cl}(R)\backslash\overline{\mathrm{Cl}}(R)\right)$ and overorders of $R$: a Bass order}{A formula for ideal classes and overorders of a Bass order}}\label{sec:formulabass}

This section proves the formula for $\#\left(\mathrm{Cl}(R)\backslash\overline{\mathrm{Cl}}(R)\right)$ and the conductor parametrization of the overorders of a Bass order $R$; see Theorem \ref{thm:bassoverorders}.

\subsection{Characterization of a Bass order}
We first extend the characterization of Bass orders in Section \ref{subsec:charbass} to the setting of Definition A.

\begin{definition}\label{def:globalbass}
Let $Z$ be a Dedekind domain with field of fractions $Q$.
Let $K$ be an \'etale $Q$-algebra.
Note that $Z$ is always $\mfo$ or $\mfo_v$ in this paper.
\begin{enumerate}
\item An order $\mathcal O$ is Gorenstein if and only if $\mathrm{cl}(\mathcal O)=\mathrm{Cl}(\mathcal O)$ by \cite[Proposition~3.4]{Ma24}.

   \item(\cite[Proposition 4.6]{Ma24} or \cite[Theorem 2.1]{LW}) An order $\mathcal{O}$ is called a Bass order if every overorder of $\mathcal{O}$ is Gorenstein, equivalently if every ideal of $\mathcal{O}$ is generated by  two elements.
    \end{enumerate}
\end{definition}

\begin{remark}\label{rmk:globalbass}
This remark is a generalization of Remark \ref{rmk:bassideal1}.
Let $\Mfo$ be an order of $K$ in the setting of Definition A and let $\Mfo'$ be an overorder of $\Mfo$.
We refer to Definition \ref{def:globalclcl} for  $\mathrm{cl}(\mathcal{O}')$ and $\overline{\mathrm{cl}(\mathcal{O}')}$.
    \begin{enumerate}
        \item In the case $Z=\mfo_v$ so that $\Mfo'$ is complete,
         \cite[Lemma 2.17]{Ma24} states that $I$ is an invertible $\mathcal{O}'$-ideal if and only if $I\otimes_{\mathcal{O}'} O'_w$ is a principal fractional $O'_w$-ideal for each $w\in |\mathcal{O}'|$.
         By Remark B.(1),
         this is equivalent to stating that  $I$ is a principal  fractional $\mathcal{O}'$-ideal. Thus
        \[     \textit{$\mathcal{O}$ is a Bass order if and only if $\#\mathrm{cl}(\Mfo')=1$ for all overorders $\Mfo'$ of $\Mfo$.}
\]
        \item In the case $Z=\mfo$, an overorder $\mathcal{O}'$  of $\Mfo$ is Gorenstein if and only if $\mathrm{cl}(\mathcal{O}')=\mathrm{Cl}(\mathcal{O}')$.
        This is equivalent to stating that $\#\overline{\mathrm{cl}(\mathcal{O}')}=1$ since
the map         $\mathrm{Cl}(\mathcal{O}) \rightarrow \mathrm{Cl}(\mathcal{O}'), I\mapsto I\Mfo',$ is surjective by \cite[Remark 3.8]{Ma20}.
        We then have the following description of a Bass order:
\[   \textit{$\Mfo$ is a Bass order if and only if $\#\overline{\mathrm{cl}(\mathcal{O}')}=1$ for all overorders $\Mfo'$ of $\Mfo$}.\]

        \item The maximal order $\mathcal{O}_K$ is a Bass order since $\overline{\mathrm{Cl}}(\mathcal{O}_K)=\mathrm{Cl}(\mathcal{O}_K)$.
\item Every order in a quadratic field extension $K/Q$ is Bass by \cite[Section~2.3]{LW}.
        \item If $\mathcal{O}$ is a Bass order, then any overorder $\mathcal{O}'$ of $\mathcal{O}$ is a Bass order as well.
\item If $Z=\mathbb{Z}$, then an order whose discriminant is fourth-power-free in $\mathbb{Z}$ is a Bass order by \cite[Theorem 3.6]{Gre82} (cf. \cite[Section 2.3]{LW}).

\item For $a\in\mathbb Z$, the cubic order
$\mathbb Z[x]/(x^3-ax^2+(a-1)x-1)$ is Bass. Its ideal class monoid
classifies Cappell--Shaneson matrices up to similarity; these matrices
are used to construct Cappell--Shaneson $4$-spheres
\cite[p.~44]{AR84}.
    \end{enumerate}
\end{remark}

\begin{proposition}\label{cor:ratioofcl}
Suppose that $R$ is a Bass order in a number field $E$. Then
\begin{enumerate}
    \item $\overline{\mathrm{Cl}}(R)=\bigsqcup\limits_{R\subset \Mfo \subset \Mfo_E}\cl(\Mfo)$ and  $\#(\mathrm{Cl}(R)\backslash\overline{\mathrm{Cl}}(R))=\textit{the number of overorders of $R$}$;

    \item $R_w$ is  a Bass order and a reduced local ring for each $w\in |R|$.
\end{enumerate}
\end{proposition}

\begin{proof}
The first claim is a restatement of Proposition \ref{prop:stratforglobal}  by Remark \ref{rmk:globalbass}(2).

For (2),    Remark \ref{rmk:globalbass}(2) yields that  $\#\overline{\mathrm{cl}(\mathcal{O})}=1$ for every overorder $\mathcal{O}$ of $R$.
   Then $\#\mathrm{cl}(\mathcal{O}\otimes_{R}R_w)=1$ by  Corollary \ref{cor:prop_global_local}(2).
   On the other hand,
  Corollary \ref{cor:prop_global_local}(1) yields that
    every overorder of $R_w$, for each $w\in |R|$, is of the form $\mathcal{O}\otimes_R R_w$ for an overorder $\mathcal{O}$ of $R$.
     Therefore $R_w$ is a Bass order  by Remark \ref{rmk:globalbass}(1).
   Reducedness follows from \cite[Lemma 2.16]{Ma24}.
\end{proof}

\begin{definition}\label{def:irreddecomsplit}
    Define two subsets $|R|^{irred}$ and $|R|^{split}$ of $|R|$  for a Bass order $R$ as follows:
\[
      \left\{
      \begin{array}{l}
|R|^{irred}\subset  \{w\in |R| : \textit{ $R_w$ is an integral domain}\};\\
|R|^{split}\subset  \{w\in |R| : \textit{ $R_w$ is not an integral domain}\}
      \end{array} \right. \textit{so that } ~~~~~~~     |R|=|R|^{irred}\sqcup |R|^{split}.
\]
\end{definition}

\begin{remark}
Henceforth, $R$ denotes a Bass order in a number field $E$. By
Proposition \ref{cor:ratioofcl}(1) and Corollary
\ref{cor:prop_global_local}(1), computing
$\#(\mathrm{Cl}(R)\backslash\overline{\mathrm{Cl}}(R))$ reduces to
counting overorders of $R_w$ for $w\in|R|$. Section
\ref{section:orbbass} treats the domain case, while the next subsection
treats the non-domain case.
\end{remark}

\subsection{Formula for orbital integrals and ideal class monoids in a Bass order: split case}\label{subsection:formulaforRsplit}

This subsection computes $\#(\Lambda_{E_w}\backslash X_{R_w})$ and $\#\overline{\mathrm{Cl}}(R_w)$ for $w\in|R|^{split}$; see Theorem \ref{prop:basssplitorbit} and Corollary \ref{cor:sorclbargeneralsplit}.
We first prove that a local Bass order $R_w$ (cf. Proposition \ref{cor:ratioofcl}(2)) contains a Bass order which is a simple extension of a DVR.

\begin{lemma}\label{lem:basssplitsimple}
For $w\in|R|^{split}$, $R_w$ contains a Bass order of $E_w$ that is a simple extension of $\Mfo_{K_w}$.
\end{lemma}

\begin{proof}
By \cite[Proposition~1.2]{Gre82}, there are discrete valuation rings
$D_1,D_2$ and an injection $\iota:R_w\hookrightarrow D_1\times D_2$
such that $\iota_i:=\operatorname{pr}_i\circ\iota:
R_w\twoheadrightarrow D_i$ for $i=1,2$.
Put $P_i=\ker(\iota_i)$. Since
$P_1P_2\subset P_1\cap P_2=\ker(\iota)=0$, every prime of $R_w$
contains one of the $P_i$. As $R_w$ is not a domain, neither $P_i$ is
zero; hence these are precisely the two minimal primes.
Thus \cite[proof of Lemma~2.10]{Ma24} gives
$E_w\cong E_1\times E_2$, where $E_i=\operatorname{Frac}(D_i)$.

  Since $R_w$ is local, each $\iota_i$ is
local and identifies the residue field of $D_i$ with that of $R_w$.
Put $\mathfrak m=\mathfrak m_{R_w}$ and
$\mathfrak n_i=\mathfrak m_{D_i}$. Then
$\iota_i(\mathfrak m)=\mathfrak n_i$ and
$\iota_i(\mathfrak m^2)=\mathfrak n_i^2$, giving surjections
$\overline\iota_i:\mathfrak m/\mathfrak m^2\to
\mathfrak n_i/\mathfrak n_i^2$.
The surjections $\overline\iota_i$ have proper kernels, and two proper
subspaces cannot cover $\mathfrak m/\mathfrak m^2$. Choose
$\bar z\notin\ker(\overline\iota_1)\cup\ker(\overline\iota_2)$ and lift
it to $z\in\mathfrak m$. Writing $\iota(z)=(a,b)$, we have
$a\in\mathfrak n_1\setminus\mathfrak n_1^2$ and
$b\in\mathfrak n_2\setminus\mathfrak n_2^2$; hence $a$ and $b$ are
uniformizers of $D_1$ and $D_2$.

The extensions $E_i/K_w$ are finite and totally ramified because their
residue fields agree. Hence every uniformizer of $D_i$ has Eisenstein
minimal polynomial over $K_w$ and generates $D_i$ over $\Mfo_{K_w}$
\cite[Ch.~I, \S6, Prop.~18]{Ser}. Let $g_1(x)$ be the minimal polynomial
of $a$ over $K_w$.

Let $\mathfrak f=(R_w:D_1\times D_2)$ be the conductor. Since
$(D_1\times D_2)/R_w$ is a finite torsion $\Mfo_{K_w}$-module, choose
$N\ge2$ such that
$\mathfrak n_1^N\times\mathfrak n_2^N\subset\mathfrak f$. For
$\delta\in\mathfrak n_2^N$, the element
$z_\delta=(a,b+\delta)=z+(0,\delta)$ lies in $R_w$, and $b+\delta$
remains a uniformizer of $D_2$. Since $g_1$ has only finitely many roots
in $E_2$ whereas $b+\mathfrak n_2^N$ is infinite, choose $\delta$ with
$g_1(b+\delta)\ne0$. If $g_2(x)$ is the minimal polynomial of
$b+\delta$ over $K_w$, then $g_1$ and $g_2$ are distinct Eisenstein
polynomials, $D_1=\Mfo_{K_w}[a]$, and
$D_2=\Mfo_{K_w}[b+\delta]$.

The assignment $x\mapsto(a,b+\delta)$ induces a homomorphism
$A:=\Mfo_{K_w}[x]/(g_1g_2)\to D_1\times D_2$.
After tensoring with $K_w$, this is an isomorphism onto
$E_1\times E_2$ by the Chinese remainder theorem. Since $A$ is
$\Mfo_{K_w}$-free, the map is injective, with image
$\Mfo_{K_w}[z_\delta]\subset R_w$.

It remains to show that $A$ is Bass. The injectivity above implies
$(g_1)\cap(g_2)=(g_1g_2)$; hence the standard exact sequence for the
ideals $(g_1)$ and $(g_2)$ gives
$0\longrightarrow A\longrightarrow D_1\times D_2
\longrightarrow \Mfo_{K_w}[x]/(g_1,g_2)\longrightarrow0$.
The product $D_1\times D_2$ is the normalization of $A$, while the
last term is a quotient of $A$ and hence cyclic as an $A$-module.
Hence Proposition
\ref{prop:intro}(4) shows that $A$ is Bass.
\end{proof}

\begin{remark}\label{rmk:rirredsplit}
Let $\mfp_1,\ldots,\mfp_r$ be the primes of $\Mfo_E$ above $w\in|R|$.
Applying
\cite[\href{https://stacks.math.columbia.edu/tag/07N9}{Lemma~07N9}]
{stacks-project} and passing to total quotient rings gives
$E_w\cong\prod_{i=1}^rE_{\mfp_i}$, where $E_{\mfp_i}$ is the
completion of $E$ at $\mfp_i$. Thus $r=1$ for
$w\in|R|^{irred}$, while the proof of Lemma
\ref{lem:basssplitsimple} gives $r=2$ for $w\in|R|^{split}$.
\end{remark}

We first treat simple split Bass orders and then pass to arbitrary ones
via Lemma \ref{lem:basssplitsimple}.

\begin{proposition}\label{lem:forbassprod}
Suppose that $R_w\cong \Mfo_{K_w}[x]/(\phi_w(x))$ (cf. Remark B.(3)).
    \begin{enumerate}
\item $\phi_w=\phi_{w,1}\phi_{w,2}$ for distinct monic polynomials
$\phi_{w,i}\in\Mfo_{K_w}[x]$ irreducible over $K_w$.
        \item  $\overline{\phi_w(x)}=(x-\overline{a})^{[E_w:K_w]}\in \kappa_{R_w}[x]$ for some $\overline{a}\in\kappa_{R_w}$.

        \item $\mathcal{O}_{E_{w,i}}\cong \Mfo_{K_w}[x]/(\phi_{w,i}(x))$ is a totally ramified extension of $\Mfo_{K_w}$ (cf. Remark B.(4)).
        \item
    $\mathcal{O}_{E_{w,i}}\cong\Mfo_{K_w}[x]/(\phi_{w,i}(x))$ is a discrete valuation ring with the commutative diagram
\[\xymatrix{R_w    \ar@{->}[d]^{\cong}   \ar@{^{(}->}[r]   & \Mfo_{E_w}\cong \mathcal{O}_{E_{w,1}} \times \mathcal{O}_{E_{w,2}} \ar@{->}[d]^{\cong}  \\ \iota: \Mfo_{K_w}[x]/(\phi_{w}(x))   \ar@{^{(}->}[r]   & \Mfo_{K_w}[x_1]/(\phi_{w,1}(x_1))\times \Mfo_{K_w}[x_2]/(\phi_{w,2}(x_2)), ~~~~~~ x\mapsto (x_1, x_2).}
\]
Here we recall that $\Mfo_{E_w}$ is the integral closure of $R_w$ in $E_w$.
    \end{enumerate}
\end{proposition}

\begin{proof}
The assertions follow from the embedding with surjective projections
constructed in the proof of Lemma \ref{lem:basssplitsimple}, the
generic-fiber factorization, and the equality of the residue fields.
\end{proof}

\begin{definition}[Relative Serre invariant]
\label{def:generalserreinv}
Let $Z$ be a PID with fraction field $Q$, let $K$ be an \'etale
$Q$-algebra, and let $\Mfo$ be an order of $K$. Set
$S_Q(\Mfo):=\operatorname{length}_Z(\Mfo_K/\Mfo)$. If
$Q/\mathbb Q_p$ is finite and $Z=\Mfo_Q$, set
$S_p(\Mfo):=\operatorname{length}_{\mathbb Z_p}(\Mfo_K/\Mfo)$; then
$S_p(\Mfo)=[\kappa_Q:\mathbb F_p]S_Q(\Mfo)$.
\end{definition}

By Proposition \ref{lem:forbassprod}(2), after translating $x$ by a
suitable element of $\Mfo_{K_w}$, we may assume that $\phi_{w,1}$ and
$\phi_{w,2}$ are Eisenstein. Put $r_i=\deg\phi_{w,i}$ and
$r=r_1+r_2$, relabel the factors so that $r_1\ge r_2$, and set
$s=\ord_{E_{w,1}}(\phi_{w,2}(x_1))$.

\begin{lemma}\label{lem:basisofsplitbass}
Under the embedding $\iota$ in Proposition
\ref{lem:forbassprod}(4), the image of $R_w$ has the
$\Mfo_{K_w}$-basis
\[
\{(x_1^j,x_2^j)\mid 0\leq j<r_2\}
\cup\{(x_1^{s+j},0)\mid 0\leq j<r_1\}.
\]
Moreover, $S_{K_w}(R_w)=s$ and
$\mathfrak f(R_w)=x_1^s\Mfo_{E_{w,1}}\times
x_2^s\Mfo_{E_{w,2}}$.
\end{lemma}

\begin{proof}
Put $s_2=\ord_{E_{w,2}}(\phi_{w,1}(x_2))$. The resultant identities give
\[
\operatorname{Res}(\phi_{w,1},\phi_{w,2})
=N_{E_{w,1}/K_w}\bigl(\phi_{w,2}(x_1)\bigr)
=(-1)^{r_1r_2}N_{E_{w,2}/K_w}
 \bigl(\phi_{w,1}(x_2)\bigr).
\]
Since both extensions $E_{w,i}/K_w$ are totally ramified, taking
normalized valuations gives $s_2=s$.

Since $\phi_{w,2}$ is monic, the classes of $x^j$ for
$0\leq j<r_2$ and $x^j\phi_{w,2}(x)$ for $0\leq j<r_1$ form an
$\Mfo_{K_w}$-basis of
$\Mfo_{K_w}[x]/(\phi_{w,1}\phi_{w,2})$. Under $\iota$, they map
respectively to $(x_1^j,x_2^j)$ and
$(x_1^j\phi_{w,2}(x_1),0)$. Since
$\phi_{w,2}(x_1)=ux_1^s$ for some
$u\in\Mfo_{E_{w,1}}^\times$, the latter vectors form an
$\Mfo_{K_w}$-basis of
$x_1^s\Mfo_{E_{w,1}}\times\{0\}$; replacing them by
$(x_1^{s+j},0)$ for $0\leq j<r_1$ gives the stated basis.

Using the injectivity of $\iota$, the standard exact sequence for
$(\phi_{w,1})$ and $(\phi_{w,2})$ identifies
$(\Mfo_{E_{w,1}}\times\Mfo_{E_{w,2}})/R_w$ with
$\Mfo_{K_w}[x]/(\phi_{w,1},\phi_{w,2})
\cong\Mfo_{E_{w,1}}/x_1^s\Mfo_{E_{w,1}}$.
The latter has $\Mfo_{K_w}$-length $s$, so $S_{K_w}(R_w)=s$.
The natural quotient maps
$\Mfo_{E_{w,i}}
\twoheadrightarrow
\Mfo_{K_w}[x]/(\phi_{w,1},\phi_{w,2})$ for $i=1,2$ have kernels
$x_1^s\Mfo_{E_{w,1}}$ and $x_2^s\Mfo_{E_{w,2}}$, respectively.
Since $R_w$ is their fiber product, multiplication by the two
idempotents of
$\Mfo_{E_{w,1}}\times\Mfo_{E_{w,2}}$ shows that  its conductor is the product of
these kernels, namely
$x_1^s\Mfo_{E_{w,1}}\times x_2^s\Mfo_{E_{w,2}}$.
\end{proof}

\begin{lemma}\label{lem:splitunitindex}
Let $\Mfo\subsetneq\Mfo_{E_w}$ be a local $\Mfo_{K_w}$-order. Then
$\#(\Mfo_{E_w}^\times/\Mfo^\times)
=(q_{R_w}-1)q_{R_w}^{S_{K_w}(\Mfo)-1}$.
\end{lemma}

\begin{proof}
Put $\kappa=\kappa_{R_w}$ and
$\mathfrak r=\mathfrak m_{E_{w,1}}\times\mathfrak m_{E_{w,2}}$.
By Proposition \ref{lem:forbassprod}(3), both factors have residue
field $\kappa$. Since $\Mfo$ contains $\Mfo_{K_w}$, the two reduction
maps $\Mfo\to\kappa$ are surjective; as $\Mfo$ is local, both have
kernel $\mathfrak m_{\Mfo}$. Hence
$\mathfrak m_{\Mfo}=\Mfo\cap\mathfrak r$ and the image of $\Mfo$ in
$\kappa\times\kappa$ is $\Delta(\kappa)$. Reduction gives exact
sequences
\[
\begin{aligned}
0&\longrightarrow \mathfrak r/\mathfrak m_{\Mfo}
 \longrightarrow \Mfo_{E_w}/\Mfo
 \longrightarrow(\kappa\times\kappa)/\Delta(\kappa)
 \longrightarrow0,\\
1&\longrightarrow(1+\mathfrak r)/(1+\mathfrak m_{\Mfo})
 \longrightarrow \Mfo_{E_w}^\times/\Mfo^\times
 \longrightarrow(\kappa^\times\times\kappa^\times)/
 \Delta(\kappa^\times)\longrightarrow1.
\end{aligned}
\]
The right-hand quotients have cardinalities $q_{R_w}$ and
$q_{R_w}-1$, respectively; by the definition of $S_{K_w}$, the first
sequence gives
$\#(\mathfrak r/\mathfrak m_{\Mfo})
=q_{R_w}^{S_{K_w}(\Mfo)-1}$.
 Let
$\mathfrak c=(\Mfo:\Mfo_{E_w})$. Since
$\Mfo\ne\Mfo_{E_w}$, one has
$\mathfrak c\subset\mathfrak m_{\Mfo}\subset\mathfrak r$.
As in Proposition \ref{rmk:counting}, the maps $1+x\mapsto x$ modulo
$\mathfrak c$ give
$\#((1+\mathfrak r)/(1+\mathfrak m_{\Mfo}))
=\#(\mathfrak r/\mathfrak m_{\Mfo})$.
The second sequence now gives the formula.
\end{proof}

\begin{theorem}\label{prop:basssplitorbit}
For $R_w\cong\Mfo_{K_w}[x]/(\phi_w(x))$ with
$w\in|R|^{split}$, we have:
\begin{enumerate}
\item $\#(\Lambda_{E_w}\backslash X_{R_w})
=q_{R_w}^{S_{K_w}(R_w)}=q_v^{S_{F_v}(R_w)}$, where
$\Lambda_{E_w}=x_1^{\mathbb Z}\times x_2^{\mathbb Z}$.
\item For $0\leq k\leq s=S_{K_w}(R_w)$, put
$\Mfo_k:=\Mfo_{K_w}[(x_1,x_2),(x_1^{s-k},0)]$.
\begin{enumerate}
\item The $\Mfo_k$ are precisely the overorders of $R_w$, and
\[
S_{K_w}(\Mfo_k)=s-k,\qquad
\mathfrak f(\Mfo_k)=x_1^{s-k}\Mfo_{E_{w,1}}
\times x_2^{s-k}\Mfo_{E_{w,2}}.
\]
\item $\#\overline{\mathrm{Cl}}(R_w)=S_{K_w}(R_w)+1
=\frac{S_{F_v}(R_w)}{[K_w:F_v]}+1$.
\end{enumerate}
\end{enumerate}
\end{theorem}

\begin{proof}
For~(1), Yun's formula \cite[Corollary~4.10]{Yun13}, Proposition
\ref{lem:forbassprod}(4), and Lemma \ref{lem:basisofsplitbass} give the
first equality; the second follows from the definitions, since
$K_w/F_v$ is unramified.

For~(2), put $Q=q_{R_w}$. The basis in Lemma
\ref{lem:basisofsplitbass} gives the displayed Serre invariant and
conductor for every $\Mfo_k$; hence they are pairwise
distinct. For $k<s$, the same basis shows that $\Mfo_k$ is local, so
Lemma \ref{lem:splitunitindex} gives
$\#(\Mfo_{E_w}^\times/\Mfo_k^\times)=(Q-1)Q^{s-k-1}$, while
$\Mfo_s=\Mfo_{E_w}$ contributes $1$.

Yun's identity \cite[Section~2.1 and p.~408]{Yun13} applies to the
finite \'etale algebra $E_w$. Grouping the ideal classes by their
multiplicator rings, as in the proof of Proposition
\ref{prop:orderidealcounting}, and using
$\operatorname{Aut}(I)=(I:I)^\times$ and
$\#\mathrm{cl}(\Mfo)=1$ from Remark \ref{rmk:globalbass}(1), gives the
first identity below; the preceding weights give the second:
\[
\begin{aligned}
\#(\Lambda_{E_w}\backslash X_{R_w})
 &=\sum_{R_w\subset\Mfo\subset\Mfo_{E_w}}
   \#(\Mfo_{E_w}^\times/\Mfo^\times),\\
\sum_{k=0}^{s}\#(\Mfo_{E_w}^\times/\Mfo_k^\times)
 &=1+(Q-1)(1+Q+\cdots+Q^{s-1})=Q^s.
\end{aligned}
\]
The second identity equals the total mass in~(1); since every summand
in the first is positive, the $\Mfo_k$ exhaust the overorders of
$R_w$. Hence $\#\overline{\mathrm{Cl}}(R_w)=s+1$ by Proposition
\ref{prop:stratforglobal} and Remark \ref{rmk:globalbass}(1).
\end{proof}

We now remove the simple-extension hypothesis.

\begin{corollary}\label{cor:sorclbargeneralsplit}
For a local Bass order $R_w$ with $w\in|R|^{split}$, one has
\[
\#(\Lambda_{E_w}\backslash X_{R_w})=q_v^{S_{F_v}(R_w)},\qquad
\#\overline{\mathrm{Cl}}(R_w)=S_{K_w}(R_w)+1.
\]
Moreover,
$\mathfrak f(R_w)=\mathfrak m_{E_{w,1}}^{S_{K_w}(R_w)}
\times\mathfrak m_{E_{w,2}}^{S_{K_w}(R_w)}$, and every overorder
$R_w'$ of $R_w$ satisfies
$R_w'=\langle R_w,\mathfrak f(R_w')\rangle$.
\end{corollary}

\begin{proof}
By Lemma \ref{lem:basssplitsimple}, choose a simple split Bass order
$R_{w,0}\subset R_w$. Write the chain in Theorem
\ref{prop:basssplitorbit} as
$R_{w,0}=\Mfo_0\subset\cdots\subset\Mfo_{s_0}=\Mfo_{E_w}$, where
$s_0=S_{K_w}(R_{w,0})$. Then $R_w=\Mfo_{k_0}$ for a unique $k_0$,
and its overorders are
$\Mfo_{k_0},\ldots,\Mfo_{s_0}$. Put
$h=s_0-k_0=S_{K_w}(R_w)$ and $Q=q_{R_w}$. Theorem
\ref{prop:basssplitorbit}(2)(a) gives the stated conductor, and there
are $h+1$ overorders; hence Proposition \ref{prop:stratforglobal} and
Remark \ref{rmk:globalbass}(1) give
$\#\overline{\mathrm{Cl}}(R_w)=h+1$.

For $k_0\leq k<s_0$, Lemma \ref{lem:splitunitindex} gives
$\#(\Mfo_{E_w}^{\times}/\Mfo_k^{\times})
=(Q-1)Q^{s_0-k-1}$, while $\Mfo_{s_0}=\Mfo_{E_w}$ contributes $1$.
The same mass identity as in the proof of Theorem
\ref{prop:basssplitorbit} therefore gives
\[
\#(\Lambda_{E_w}\backslash X_{R_w})
=1+(Q-1)(1+Q+\cdots+Q^{h-1})
=Q^h=q_v^{S_{F_v}(R_w)}.
\]
Finally, if $R_w'=\Mfo_k$, then $\mathfrak f(R_w')$ contains
$(x_1^{s_0-k},0)$, which generates $\Mfo_k$ over $R_{w,0}$. Thus
$\Mfo_k\subset\langle R_w,\mathfrak f(R_w')\rangle$; the reverse
inclusion follows from $R_w,\mathfrak f(R_w')\subset\Mfo_k$, proving
the final assertion.
\end{proof}

\subsection{\texorpdfstring{Formula for $\#(\mathrm{Cl}(R)\backslash\overline{\mathrm{Cl}}(R))$}{Formula for ideal classes modulo invertible ideals}}
Combining Sections \ref{subsec:1stcase}, \ref{subsec:2ndcase}, and
\ref{subsection:formulaforRsplit}, we obtain the global conductor
parametrization of the overorders of $R$ and the formula for
$\#(\mathrm{Cl}(R)\backslash\overline{\mathrm{Cl}}(R))$.

Since the conductor is the annihilator of the $R$-module $\Mfo_E/R$,
\cite[\href{https://stacks.math.columbia.edu/tag/07T8}{Lemma 07T8}]
{stacks-project} gives
$\mathfrak f(R)\otimes_RR_w=\mathfrak f(R_w)$ for $w\in|R|$.
Thus Remark \ref{rmk:rirredsplit}, \eqref{eq:relsercon}, and Corollary
\ref{cor:sorclbargeneralsplit} give
\begin{equation}\label{descofconductorofR}
\mathfrak f(R)=
\left(\mathfrak p_1^{2l_1}\cdots\mathfrak p_r^{2l_r}\right)
\left(\mathfrak q_1^{m_1}\cdots\mathfrak q_s^{m_s}\right)
\left((\mathfrak r_1\widetilde{\mathfrak r}_1)^{n_1}\cdots
(\mathfrak r_t\widetilde{\mathfrak r}_t)^{n_t}\right).
\end{equation}
Here the displayed prime ideals of $\Mfo_E$ are pairwise distinct:
$\mathfrak p_i$ lies over $w_i\in|R|^{irred}$ with
$[\kappa_{E_{w_i}}:\kappa_{R_{w_i}}]=1$; $\mathfrak q_j$ lies over
$w_j\in|R|^{irred}$ with
$[\kappa_{E_{w_j}}:\kappa_{R_{w_j}}]=2$; and
$\mathfrak r_k,\widetilde{\mathfrak r}_k$ lie over
$w_k\in|R|^{split}$. The exponents satisfy
\begin{equation}\label{eq:relabetweenlocalsercon}
l_i=S_{K_{w_i}}(R_{w_i}),\qquad
m_j=S_{K_{w_j}}(R_{w_j}),\qquad
n_k=S_{K_{w_k}}(R_{w_k}).
\end{equation}

\begin{theorem}\label{thm:bassoverorders}
Let $R$ be a Bass order in a number field $E$. Then:
\begin{enumerate}
\item \(\displaystyle
\begin{aligned}[t]
\#\bigl(\mathrm{Cl}(R)\backslash\overline{\mathrm{Cl}}(R)\bigr)
&=\prod_{w\mid\mathfrak f(R)}\bigl(S_{K_w}(R_w)+1\bigr)\\
&=\prod_{p\in\mathcal P_R}\prod_{w\mid p}
\left(\frac{S_p(R_w)}{d_{R_w}}+1\right)
=\prod_{i=1}^{r}(l_i+1)\prod_{j=1}^{s}(m_j+1)
\prod_{k=1}^{t}(n_k+1).
\end{aligned}\)\par
Here $d_{R_w}=[\kappa_{R_w}:\mathbb F_p]$ and
$\mathcal P_R$ is the set of rational primes dividing
$\operatorname{disc}(R)/\operatorname{disc}(\Mfo_E)$.

\item Every overorder $R'$ of $R$ is of the form
$\langle R,I_{l_i',m_j',n_k'}\rangle$ for a unique ideal
\[
\begin{gathered}
I_{l_i',m_j',n_k'}=
\left(\mathfrak p_1^{2l_1'}\cdots\mathfrak p_r^{2l_r'}\right)
\left(\mathfrak q_1^{m_1'}\cdots\mathfrak q_s^{m_s'}\right)
\left((\mathfrak r_1\widetilde{\mathfrak r}_1)^{n_1'}\cdots
(\mathfrak r_t\widetilde{\mathfrak r}_t)^{n_t'}\right),\\
0\le l_i'\le l_i,\qquad 0\le m_j'\le m_j,\qquad
0\le n_k'\le n_k.
\end{gathered}
\]
In this case $I_{l_i',m_j',n_k'}=\mathfrak f(R')$, and hence
$R'=\langle R,\mathfrak f(R')\rangle$.
\end{enumerate}
\end{theorem}

\begin{proof}
\begin{enumerate}
\item Proposition \ref{prop:global_local}, Theorem
\ref{thm:alllocalresults}, and Corollary
\ref{cor:sorclbargeneralsplit} give the first product in~(1), since
its nontrivial factors are precisely those with
$w\mid\mathfrak f(R)$. The relation
$S_p(R_w)=d_{R_w}S_{K_w}(R_w)$ gives the second product after grouping
by the rational prime below $w$. The standard index--discriminant
formula identifies these primes with $\mathcal P_R$, and
\eqref{eq:relabetweenlocalsercon} gives the last product.

\item For each admissible tuple, put $I=I_{l_i',m_j',n_k'}$ and
$R_I=\langle R,I\rangle$. Since $(R_I)_w=\langle R_w,I_w\rangle$,
Theorem \ref{thm:alllocalresults} and Corollary
\ref{cor:sorclbargeneralsplit} give
$\mathfrak f((R_I)_w)=I_w$. Compatibility of conductors with completion
gives $\mathfrak f(R_I)=I$, so the $R_I$ are pairwise distinct. Their
number is the last product in~(1), hence Proposition
\ref{cor:ratioofcl}(1) shows that they exhaust all overorders. \qedhere
\end{enumerate}
\end{proof}

\begin{remark}\label{rmk:reinterpmainresult}
For a Bass order $R$, the standard index--discriminant formula gives
the first equality below. Since $R$ is Gorenstein,
\cite[Corollary~4, p.~84]{DCD} gives the second:
\[
\#(\Mfo_E/R)^2
=\left|\frac{\operatorname{disc}(R)}
{\operatorname{disc}(\Mfo_E)}\right|
=N_{E/\mathbb Q}(\mathfrak f(R)).
\]
Moreover, $\Mfo_E/R\cong
\prod_{w\mid\mathfrak f(R)}\Mfo_{E_w}/R_w$, so the definition of
$S_{K_w}$ gives
\[
\#(\Mfo_E/R)=
\prod_{w\mid\mathfrak f(R)}q_{R_w}^{S_{K_w}(R_w)},\qquad
N_{E/\mathbb Q}(\mathfrak f(R))=
\prod_{w\mid\mathfrak f(R)}q_{R_w}^{2S_{K_w}(R_w)}.
\]
These identities agree with \eqref{eq:relabetweenlocalsercon}.
\end{remark}

\end{document}